\documentclass{article}
\usepackage{amsmath}
\usepackage{amsthm,amssymb}
\usepackage{epsfig,graphics}
\usepackage{psfrag}
 \columnwidth\textwidth

\newcommand{\smfrac}[2]{{\frac{#1}{#2}}}
\def\eps{\varepsilon}
\def\w{\Omega}
\def\what{\widehat{\Omega}}
\def\Vhat{\widehat{V}}
\def\Vt{\check{V}}
\def\Va{V(q_1,0)}

\def\Vc{V(q_1 + h P_1, 0)}
\def\Vd{V(q_1,\eps (x \sin \sigma - Y \cos \sigma))}
\def\Ve{V(q_1 + h P_1,\eps (x \sin \Sigma - Y \cos \Sigma))}

\def\Xint#1{\mathchoice
	{\XXint\displaystyle\textstyle{#1}}%
	{\XXint\textstyle\scriptstyle{#1}}%
	{\XXint\scriptstyle\scriptscriptstyle{#1}}%
	{\XXint\scriptscriptstyle\scriptscriptstyle{#1}}%
\!\int}
\def\XXint#1#2#3{{\setbox0=\hbox{$#1{#2#3}{\int}$ }
	\vcenter{\hbox{$#2#3$ }}\kern-.6\wd0}}

\def\dashint{\Xint-}
\def\SFD{S_{\eps, {\rm FD}}}
\def\Hsimp{H_2}
\def\H3simp{H_3}

\newcommand{\RR}{\mathbb R}
\newcommand{\ZZ}{\mathbb Z}
\newcommand{\NN}{\mathbb N}

\theoremstyle{plain}
\newtheorem{lemma}{Lemma}[section]
\theoremstyle{remark}
\newtheorem{remark}{Remark}[section]

\newcommand{\dps}{\displaystyle}
\newcommand{\Id}{\mbox{Id}}
\newcommand{\alp}{r}

\newlength{\minipagewidth}
\newcommand{\bookbox}[1]{
\par\medskip\noindent
\framebox[\textwidth]{
\begin{minipage}{\minipagewidth}
{#1}
\end{minipage} } \par\medskip }

\newtheorem{algorithm}{Algorithm}

\begin{document}

\title{Symplectic schemes for highly oscillatory Hamiltonian systems:
  the homogenization approach\\ beyond the constant frequency case}
\author{Matthew Dobson, Claude Le Bris, and Fr\'ed\'eric Legoll}
\date{\today}

\maketitle

\begin{abstract}
We follow up on our previous works~\cite{cras1,dcds} which presented a possible
approach for deriving symplectic schemes for a certain class 
of highly oscillatory
Hamiltonian systems. The approach considers the
Hamilton-Jacobi form of the equations of motion, formally homogenizes it
and infers an appropriate symplectic integrator for the original system. 
In~\cite{cras1,dcds} the case of a system exhibiting a single
constant fast frequency was considered. The present work successfully
extends the approach to systems that have either one \emph{varying} fast
frequency or \emph{several constant} 
frequencies. Some related issues are also examined. 
\end{abstract}

\tableofcontents

\section{Introduction}

The purpose of this work is to describe an approach for 
constructing symplectic numerical
integrators for systems within the following class of highly
oscillatory Hamiltonian systems:
\begin{equation}
\label{edo}
\frac{dq}{dt} = \frac{\partial H_\eps}{\partial p}(q,p),
\quad 
\frac{dp}{dt} = -\frac{\partial H_\eps}{\partial q}(q,p),
\end{equation}
with
\begin{equation}
\label{ham_general}
H_\eps(q_1,q_2,p_1,p_2) = 
\frac{p_1^T p_1}{2} + \frac{p_2^T p_2}{2} +
V(q_1,q_2) + \frac{q_2^T \Omega(q_1)^2 q_2}{2\eps^2} ,
\end{equation}
where $\eps \ll 1$ is a small parameter, where $q=(q_1,q_2)\in
\RR^s\times\RR^f$ and $p=(p_1,p_2)\in
\RR^s\times\RR^f$ (the upper indices $s$ and $f$ designate the slow
and fast variables respectively), and where the interaction potential
$V$ and the fast frequency factor $\Omega(q_1)$ are all
independent of $\eps.$  We assume that $\Omega(q_1)$ is a symmetric matrix
whose eigenvalues are positive and bounded away from zero and that $V$
is bounded from below. We also assume that the initial conditions for
(\ref{edo}) scale with $\eps$ so that the initial energy is bounded
independently of $\eps$.  
Since the small parameter $\eps$ imparts fast frequencies on the model,
direct numerical simulation of~\eqref{edo} to times of $O(1)$ and beyond
is computationally expensive.
This type of Hamiltonian system has already been studied in many works
\cite{schuette,rennes,grimm,HLW,flavors}. See~\cite{cohen_review,dcds} for a short
review of the literature, and \cite[Chap. XIII and XIV]{HLW} for a
general overview.  This paper continues the work initiated
in~\cite{cras1,dcds} by extending the generality of examples treated
in the Hamilton-Jacobi framework.  
The major assumption remaining on the form of~\eqref{ham_general} is
that the leading order behavior in the fast variables is that of a
harmonic oscillator 
$$
\frac{p_2^T p_2}{2} + \frac{q_2^T \Omega(q_1)^2 q_2}{2\eps^2}.
$$
In Section~\ref{sec:pendulum}, we briefly discuss the use
of coordinate transforms for well-behaved, though fully non-linear,
potentials~\cite[Sec. XIV.3]{HLW}.  The remainder of that section
describes the use of our approach on a particular Hamiltonian whose fast potential 
does not take the form of a harmonic oscillator in the original
coordinates.

As in our previous works~\cite{cras1,dcds}, our approach to the problem
is based on the Hamilton-Jacobi formalism. 
Let $S_\eps(t,q,P)$ be the solution to
\begin{equation}
\label{hj}
\partial_t S_\eps = H_\eps \left( q + \partial_{P} S_\eps,P
\right), \quad 
S_\eps(0,q,P) = 0. 
\end{equation}
For all $(q,p,t)$, it is known that the functions $(Q(t),P(t))$, implicitly defined 
by 
\begin{equation}
\label{symp}
p = P(t) + \frac{\partial S_\eps}{\partial q}(t,q,P(t)), 
\quad \quad
Q(t) = q + \frac{\partial S_\eps}{\partial P}(t,q,P(t)), 
\end{equation}
are solutions to (\ref{edo}) with initial conditions $(q,p)$.  
For any approximation $\widetilde{S}_\eps$ of the generating
function and stepsize $h$, there is a corresponding symplectic map 
$\Psi_h : (q,p) \mapsto (Q(h),P(h))$ defined by the implicit relations
\begin{equation}
\label{symp_scheme}
p = P(h) + \frac{\partial \widetilde{S}_\eps}{\partial q}(h,q,P(h)), 
\quad \quad
Q(h) = q + \frac{\partial \widetilde{S}_\eps}{\partial P}(h,q,P(h)).
\end{equation}
In the following we construct a
function $\widetilde{S}_\eps$ that approximates the solution
of~\eqref{hj} for small $t$ and $\eps$.  As
long as we solve~\eqref{symp_scheme} exactly,
this results in a symplectic numerical scheme (see~\cite{feng}).
Being able to generate a class of symplectic schemes motivates the 
strategy, which we use in the sequel, of
making all {\em approximations} on the level of the generating
function $\widetilde{S}_\eps$ and from there 
solving~\eqref{symp_scheme} to build the
numerical scheme.  
Of course, we cannot expect to solve these equations exactly, 
but we do so to
high precision and observe good behavior with respect to preserving
invariants.

We work
within the parameter regime $\eps \ll h \ll \eps^{1/3}.$  This yields
a computational speed-up in comparison to standard algorithms such as
velocity Verlet, where the time step must be taken smaller than the
characteristic time scale of the fast motion, $O(\eps)$.  
However, we note that the energy preservation property of symplectic schemes is 
typically proven in the limit $h\rightarrow 0$ for a given Hamiltonian,
and we are thus working outside the theoretical regime of 
symplectic schemes.  We nonetheless choose symplecticity
as a goal for our scheme and provide numerical tests of the preservation
of energy and other invariants.

In contrast to~\cite{cras1,dcds}, where only the case of a system
exhibiting \emph{one single constant} fast frequency was considered,
the present work successfully extends the approach to systems that
have either {\em one varying} fast frequency, or {\em several
constant} frequencies.  Many of the results presented here have been
announced in~\cite{cras2}. 

In addition to extending the approach to these two more general
settings, we address one issue that we left pending in our previous
works~\cite{cras1,dcds}.  This issue is related to the approximation
of the high order derivatives present in the numerical scheme. One
drawback of  symplectic integrators constructed using the
Hamilton-Jacobi approach is that they require high order derivatives
of the Hamiltonian. The question then arises of whether or not these
high order derivatives can be approximated by the corresponding finite
differences without sacrificing the other advantages of the approach
(symplecticity, accuracy, \ldots). We demonstrate here that it is
indeed possible to use such finite difference approximations and keep
all the features of the algorithm up to a chosen order. 

\medskip

The present article is organized as follows. 
In Section~\ref{sec:approach}, we consider a Hamiltonian of the form 
\begin{equation}
\label{ham}
H_\eps(q_1,q_2,p_1,p_2) = \frac{p_1^T p_1}{2} + \frac{p_2^T p_2}{2} +
V(q_1,q_2) + \Omega(q_1)^2 \frac{q_2^T q_2}{2\eps^2},
\end{equation}
where $\w(q_1)$ is a {\em scalar} fast frequency that {\em depends on the
slow variables}.  As was already mentioned for the constant frequency
case in~\cite{dcds}, the best option to implement our approach is to
precondition the fast motion using a change of variables. This is to
be performed prior to writing the Hamilton-Jacobi equation associated
to the Hamiltonian dynamics. However, the dependency of $\Omega$ upon
the slow variables introduces substantial new difficulties in this
preconditioning step as compared to~\cite{dcds}, which will be described in
details below.  

The algorithm that we obtain following this strategy has been
introduced in~\cite{cras2}. In this paper we provide a comprehensive
discussion of its derivation and properties as well as numerical tests
demonstrating that the algorithm has favorable error performance in
terms of  resonances.  The algorithm's computational efficiency is
comparable, although slightly lower due to implicitness, to another
established approach for small $\eps$ (such as $10^{-4}$) and does not
break down as $\eps \rightarrow 0.$  We further show that the exchange
of actions (energy divided by frequency) among the fast degrees of
freedom is captured well by our algorithm.  The algorithm presented
has many possible variants based on the choice of approximation of the
generating function in~\eqref{hj}. Among the possible variants, are
those that reduce to the algorithms given in~\cite{dcds} whenever
$\Omega$ is constant. 

\medskip

In Section~\ref{sec:multi}, we next consider the case of a {\em
  non-scalar constant} fast frequency $\w,$ 
\begin{equation}
\label{ham_matrix}
H_\eps(q_1,q_2,p_1,p_2) = \frac{p_1^T p_1}{2} + \frac{p_2^T p_2}{2} +
V(q_1,q_2) + \frac{q_2^T \w^2 q_2}{2\eps^2}.
\end{equation}
We consider first the case when the fast frequencies present in $\Omega$
are non-resonant, and next address the case when some frequencies are
resonant. These terms will be defined in Section~\ref{sec:multi}. We will show the
efficiency of the algorithm obtained on a classical test case.

\medskip

We conclude our work by studying an example system composed of a single
point attached by a spring to a pivot in the plane (see Section~\ref{sec:pendulum}). 
Such a system cannot initially be written as a system with Hamiltonian of the form~\eqref{ham_general}. 
After a change of coordinates, we transform the
Hamiltonian into the form
\begin{equation}
\label{ham_mass}
H_\eps(q_1,q_2,p_1,p_2) = \frac{p_1^T M(q_2)^{-1} p_1}{2} + \frac{p_2^T p_2}{2} +
V(q_1,q_2) + \w^2 \frac{q_2^T q_2}{2\eps^2} ,
\end{equation}
with a non-constant mass matrix $M(q_2)$. Form~\eqref{ham_mass} is now
compatible with~\eqref{ham_general}, up to the fact that the mass matrix
$M(q_2)$ depends on the fast position $q_2$. We show that our general approach
again applies, yielding an algorithm with very good numerical
performance. 

\bigskip

%%%%%%%%%%%%%%%%%%%%%%%%%%%%%%%%%%%%%%%%%%%%%%%%%%%%%%%%%%%%%%%%%%%%%%%%%
\noindent{\bf Notation}
Before proceeding, we briefly fix the following notation.
For vectors $u,v \in \RR^n,$ we define the dot product 
$u \cdot v \in \RR$ by
\begin{equation*}
u \cdot v = u^T v = \sum_{i=1}^n u_i v_i,
\end{equation*}
where we also use as shorthand $u^2 = u \cdot u.$
For $u \in \RR^m$ and $v \in \RR^n,$ we define the tensor product 
$u \otimes v \in \RR^{m \times n}$ by
\begin{equation*}
(u \otimes v)_{ij} = u_i v_j, \qquad \text{ for } 1\leq i \leq m, \quad 1\leq j\leq n.
\end{equation*} 
For matrices $A,B \in \RR^{m\times n}$, we define $A:B \in \RR$ by
\begin{equation*}
A:B = \sum_{i=1}^m \sum_{j=1}^n A_{ij} B_{ij}.
\end{equation*} 
We employ two different derivative notations.  For functions
$S(q_1,P_1)$ and $V(q_1,q_2)$ with $q_1,P_1 \in \RR^s$ and $q_2 \in
\RR^f$, we denote the entries of $q_1$ by $q_1 =
(q_{1,1},\dots,q_{1,s})$, and we write  
$\partial_{q_1} S$ and $\nabla_1 V$ to denote the vectors
\begin{equation*}
(\partial_{q_1} S)_j = \frac{\partial S}{\partial q_{1,j}}, \qquad
(\nabla_{1} V(q_1,q_2))_j = \frac{\partial V}{\partial q_{1,j}}(q_1,q_2)
\qquad 1 \leq j \leq s.
\end{equation*}

%%%%%%%%%%%%%%%%%%%%%%%%%%%%%%%%%%%%%%%%%%%%%%%%%%%%%%%%%%%%%%%%%%%%%%%%%
\section{The case of a scalar non-constant frequency}
\label{sec:approach}

In this section, we consider the case~\eqref{ham} of a Hamiltonian with
a scalar 
non-constant fast frequency $\w(q_1)$. Slightly changing the notation in
comparison to~\eqref{ham}, the Hamiltonian reads
\begin{equation}
\label{ham_encore}
H_\eps(\check{q}_1,\check{q}_2,\check{p}_1,\check{p}_2) = 
\frac{\check{p}_1^T \check{p}_1}{2} + 
\frac{\check{p}_2^T \check{p}_2}{2} +
\Vt(\check{q}_1,\check{q}_2) + 
\Omega(\check{q}_1)^2 \frac{\check{q}_2^T  \check{q}_2}{2\eps^2} .
\end{equation}
We first precondition the equations with a change of variables that
takes into account the most oscillatory component of the fast
solution. Following~\cite[page 555]{HLW}, we introduce the 
symplectic change of  variables
\begin{equation}
\label{cov}
q_1= \check{q}_1, \qquad
p_1 = \check{p}_1 - \frac{\nabla_1 \w(\check{q}_1)}{2 \w(\check{q}_1)} 
\check{q}_2^T \check{p}_2, 
\qquad 
q_2 = \frac{\sqrt{\w(\check{q}_1)}}{\sqrt{\eps}} \check{q}_2, \qquad
p_2 =\frac{\sqrt{\eps}}{\sqrt{\w(\check{q}_1)}} \check{p}_2,
\end{equation}
which transforms~\eqref{ham_encore} into
$$
H_1(q_1,q_2,p_1,p_2) = 
\frac{1}{2} \left(p_1 +\frac{\nabla_1 \w(q_1) q_2^T p_2}{2 \w(q_1) } \right)^2 +
\frac{\w(q_1)}{2\eps}
(p_2^T p_2 + q_2^T q_2) + \Vt\left(q_1, \frac{\sqrt{\eps}}{\sqrt{\w(q_1)}}
q_2\right).
$$
We neglect part of the slow momentum term and rewrite the slow potential
to arrive at 
$$
\Hsimp(q_1,q_2,p_1,p_2) = \frac{1}{2} p_1^T p_1 + 
\frac{\w(q_1)}{2\eps} (p_2^T p_2 + q_2^T q_2) + 
V(q_1, \sqrt{\eps} q_2),
$$
where 
$$
V(q_1,q_2) = \Vt\left(q_1, \frac{1}{\sqrt{\w(q_1)}}
q_2\right).
$$
The system of equations corresponding to $\Hsimp$ is
\begin{equation*}
%\label{H2_dynamics}
\begin{split}
\dot{q}_1 &= p_1,
\quad \quad
\dot{p}_1 = -\nabla_1 V(q_1, \sqrt{\eps} q_2) - \frac{\nabla_1 \w(q_1)}{2 \eps} 
         (p_2^T p_2 + q_2^T q_2)
\\
\dot{q}_2 &= \frac{\w(q_1)}{\eps} p_2,
\quad \quad
\dot{p}_2 = - \sqrt{\eps} \ \nabla_2 V(q_1, \sqrt{\eps} q_2) - \frac{\w(q_1)}{\eps} q_2
\end{split}
\end{equation*}

In the spirit of~\cite{dcds}, we introduce the fast time 
\begin{equation}
\label{eq:sigma}
\sigma(t) = \frac{1}{\eps} \int_0^t \w(q_1(s)) \, ds
\end{equation} 
and consider the time-dependent change of variables $(q_2,p_2) \mapsto
(\widetilde{x},\widetilde{y})$ defined by
$$
q_2 = \widetilde{x} \cos(\sigma(t)) + \widetilde{y} \sin(\sigma(t)),
\quad \quad
p_2 = -\widetilde{x} \sin(\sigma(t)) + \widetilde{y} \cos(\sigma(t)),
$$
which captures the highest frequency component of the fast
variables.  The transformation is motivated by the fact that if
$V$ were independent of $q_2,$ then $\widetilde{x}$ and
$\widetilde{y}$ would be constant.  
This gives the dynamics
\begin{equation}
\label{turned_dynamics}
\begin{split}
\dot{q}_1 &= \frac{\partial \H3simp}{\partial p_1},
\quad \quad
\dot{p}_1 = -\frac{\partial \H3simp}{\partial q_1} -
\frac{\nabla_1 \w(q_1)}{2 \eps} 
(\widetilde{x}^T \widetilde{x}+\widetilde{y}^T \widetilde{y}),
\\
\dot{\widetilde{x}} &= \frac{\partial \H3simp}{\partial \widetilde{y}},
\quad \quad
\dot{\widetilde{y}} = -\frac{\partial \H3simp}{\partial \widetilde{x}},
\end{split}
\end{equation}
where
$$
\H3simp(q_1,\widetilde{x},p_1,\widetilde{y},t) = \frac{1}{2} p_1^2 + 
V\left[ q_1, \sqrt{\eps} (\widetilde{x} \cos \sigma(t) +
\widetilde{y} \sin \sigma(t)) \right].
$$
Note that the dynamics on $(q_1, \widetilde{x}, p_1, \widetilde{y})$ includes
a {\emph {memory}} term due to the dependence of $\sigma$ on the
history of $q_1,$ see \eqref{eq:sigma}. 
On the other hand, in the case of
a constant frequency, $\w,$ there is no memory term since $\sigma$ does not
depend on $q_1$, and~\eqref{turned_dynamics}
is the Hamiltonian dynamics associated to $\H3simp$. In the present case, 
rather than operate on a system with memory, we consider $\sigma$
as an additional variable and add the corresponding conjugate variable 
$\dps \widetilde{a} = \frac{1}{2} (\widetilde{x}^T \widetilde{x} + \widetilde{y}^T \widetilde{y}).$ 
We then have the dynamics
\begin{equation*}
\begin{split}
\dot{q}_1 &= p_1, 
\\
\dot{p}_1 &= -\nabla_1 V(q_1, \sqrt{\eps} (\widetilde{x} \cos \sigma +
\widetilde{y} \sin \sigma)) - \frac{\nabla_1 \w(q_1)}{\eps} \,
\widetilde{a} ,
\\
\dot{\widetilde{x}} &= \hphantom{-} \sqrt{\eps} \sin \sigma \ \nabla_2 V(q_1, 
\sqrt{\eps} (\widetilde{x} \cos \sigma + \widetilde{y} \sin \sigma)),
\\
\dot{\widetilde{y}} &= - \sqrt{\eps} \cos \sigma \ \nabla_2 V(q_1, 
      \sqrt{\eps} (\widetilde{x} \cos \sigma + \widetilde{y} \sin
      \sigma)),
\\
\dot \sigma &=  \frac{1}{\eps} \w(q_1), 
\\
\dot{\widetilde{a}} &= \sqrt{\eps} \
(\widetilde{x} \sin \sigma - \widetilde{y} \cos\sigma)^T
\nabla_2 
V(q_1, \sqrt{\eps} (\widetilde{x} \cos \sigma + \widetilde{y} \sin \sigma)),
\end{split}
\end{equation*}
which is a Hamiltonian dynamics associated with the energy
\begin{equation}
\label{h4}
H_4(q_1,\widetilde{x},\sigma, p_1,\widetilde{y},\widetilde{a}) = \frac{1}{2} p_1^T p_1 + V(q_1, 
          \sqrt{\eps} (\widetilde{x} \cos \sigma + \widetilde{y} \sin \sigma))
		+ \widetilde{a} \frac{\w(q_1)}{\eps}.
\end{equation}
 
\bigskip

In the sequel, we take the above Hamiltonian $H_4$ as a starting point
for our manipulations. We hence write the Hamilton-Jacobi equation
associated to $H_4$ and then rescale the variables $\widetilde{a},
\widetilde{x},$ and $\widetilde{y}$, 
$$
a = \frac{\widetilde{a}}{\eps}, \ x = \frac{\widetilde{x}}{\sqrt{\eps}}, \text{ and } 
y = \frac{\widetilde{y}}{\sqrt{\eps}} ,
$$
so that $(q_1,x,P,Y,a)$ are of the same order.  We choose not to rescale
$\Sigma$ despite the fact that it is $O(t/\eps)$ at time $t$ because
it plays the role of a fast time in the original dynamics.  The Hamilton-Jacobi equation becomes 
\begin{equation}
\begin{split} 
\partial_t S_\eps(t,q_1,x,\Sigma,P_1,Y,a) 
&
= H_4 \left[q_1 + \partial_{P_1} S_\eps, \sqrt{\eps}\left(x + \frac{1}{\eps}
\partial_Y S_\eps\right), 
\Sigma, P_1, \sqrt{\eps} Y, \eps a - \partial_\Sigma S_\eps \right]
\\
&= \frac{1}{2} P_1^2
+ V\big(q_1+\partial_{P_1} S_\eps, 
( \eps x + \partial_Y S_\eps) \cos\Sigma + \eps Y \sin\Sigma\big) 
\\
&\qquad + \left(a - \frac{1}{\eps} \partial_\Sigma S_\eps \right) 
\w(q_1+\partial_{P_1} S_\eps), 
\\
S_\eps(0,q_1,x,\Sigma,P_1,Y,a) &= 0.
\end{split} 
\label{hj2}
\end{equation}
The equations corresponding to~\eqref{symp} are
\begin{equation} 
\label{step}
\begin{split}
Q_1 &= q_1 + \partial_{P_1} S_\eps, \qquad \qquad 
X = x + \frac{1}{\eps}\partial_{Y} S_\eps, \qquad \qquad
\Sigma = \sigma + \frac{1}{\eps}\partial_{a} S_\eps,\\
P_1 &= p_1 - \partial_{q_1} S_\eps, \qquad \qquad 
Y = y - \frac{1}{\eps} \partial_{x} S_\eps, \qquad \qquad
A = a - \frac{1}{\eps}\partial_{\Sigma} S_\eps,
\end{split}
\end{equation}
where all derivatives of $S_\eps$ are evaluated at
$(t,q_1,x,\Sigma,P_1,Y,a).$

%%%%%%%%%%%%%%%%%%%%%%%%%%%%%%%%%%%%%%%%%%%%%%%%%%%%%%%%%%%%%%%%%%%%%%%%%
\subsection{Expanding in $\eps$}

We wish to expand the generating function as a series in powers of
$\eps$. Recall that there are two balancing factors in choosing an
ansatz.  Since the expanded Hamilton-Jacobi equation of order $\eps^k$
will involve the generating function up to order $\eps^{k+1},$  we
need to make sufficiently strong assumptions to be able to separate out the
$O(\eps^{k+1})$ terms from the lower order terms into independent
equations and thereby close the hierarchy. On the other hand, we also
need to have a sufficiently general form in order to satisfy the resulting
equations.  We make the two-scale ansatz
\begin{equation}
\label{ansatz}
S_\eps(t, q_1, x, \Sigma, P_1, Y, a) = S_0(t,q_1,P_1,a) +  
\sum_{k=1}^\infty \eps^k S_k\left(t, \Sigma - \frac{1}{\eps}
\tau(t,q_1,P_1,a), q_1, x, \Sigma, P_1, Y, a\right)
\end{equation} 
where we assume $S_k(t, \beta, q_1, x, \gamma, P_1, Y, a)$ is $2\pi$ 
periodic in $\gamma$ and choose the fast time $\tau$ to satisfy 
\begin{equation}
\label{tau}
\begin{split}
&\partial_t \tau = \big(\nabla_1 V(q_1+ \partial_{P_1} S_0, 0) 
+ a \nabla_1 \w(q_1+ \partial_{P_1} S_0)\big) \cdot \partial_{P_1} \tau
+ \w(q_1+ \partial_{P_1} S_0), \\
&\qquad \qquad \tau(0,q_1,P_1,a) =0.
\end{split}
\end{equation}

The periodicity assumption on $\gamma$ is needed to separate the orders
in the  
resulting hierarchy of equations and is justified both by the fact that
 equation~\eqref{hj2} is $2 \pi$-periodic in $\Sigma$ and that
the generating function of a harmonic oscillator is periodic in
time. The simplifying choice of $S_0$ independent of $\Sigma$ is
consistent with~\eqref{step} and the fact that $A$ is $O(1)$ over long
times. Similar arguments justify the choice that $S_0$ is
a function of $(t,q_1,P_1,a)$ only. Note that, as is common in
homogenization expansions, we have introduced a fast time
$\tau/\eps$ in $S_k$, which is considered to be independent from the
slow time $t$.  Because of the non-constant frequency $\w(q_1)$, we
choose for $\tau$ to depend on the slow variables $(t,q_1,P_1,a).$
What is unusual is that we introduced the fast time as part of a
transport term $\dps \Sigma - \tau/\eps$.  We motivate this by briefly
considering the Hamiltonian $H(a,\sigma) = - \eps \cos(\sigma) - a/\eps$,
a highly simplified version of $H_4$ (see~\eqref{h4}).  It is periodic
in $\sigma$ and keeps the conjugate variable $a,$ with frequency $\w=1$.  
In this case, the generating function
$S_\eps(t,a,\Sigma) = - a t/\eps - \eps^2 [\sin \Sigma - \sin(\Sigma -
t/\eps)]$ satisfies the Hamilton-Jacobi equation $\partial_t S =
H(a+\partial_\Sigma S, \Sigma)$ with initial condition $S(0,a,\Sigma)
= 0.$  Inspired by the transport structure in $\sin(\Sigma -
t/\eps),$
we choose the independent variables $\Sigma$ and
$\dps \Sigma - \tau/\eps$ for $S_k$, rather than $\Sigma$ and $\tau/\eps$. 

\smallskip

In view of~\eqref{ansatz}, we then have
\begin{equation*}
\begin{array}{rclrcl} 
\partial_t S_\eps &=& \sum_{k=0}^{\infty} \eps^k 
(\partial_t S_k - \partial_t \tau \  \partial_\beta S_{k+1}), 
&
\partial_Y S_\eps &=& \sum_{k=1}^{\infty} \eps^k 
\partial_Y S_k,
\\
\partial_{P_1} S_\eps &=& \sum_{k=0}^{\infty} \eps^k 
(\partial_{P_1} S_k - \partial_{P_1} \tau \ \partial_\beta S_{k+1}),
&
\partial_\Sigma S_\eps &=& \sum_{k=1}^{\infty} \eps^k 
(\partial_\beta S_k + \partial_\gamma S_k).
\end{array}
\end{equation*} 
From the initial condition $S_\eps(0,q_1,x,\Sigma,P_1,Y,a) = 0$ and the 
ansatz, we have
\begin{equation}
\label{ICs}
S_0(0,q_1,P_1,a) =0
\quad \text{and} \quad
S_k(0,\gamma,q_1,x,\gamma,P_1,Y,a) = 0
\ \ \text{ for } 1 \leq k < \infty,
\end{equation}
where we note the repetition of arguments in $S_k$ results
from~\eqref{ansatz} and the fact that $\tau(0,q_1,P_1,a)=0.$

%%%%%%%%%%%%%%%%%%%%%%%%%%%%%%%%%%%%%%%%%%%%%%%%%%%%%%%%%%%%%%%%%%%%%%%%%
\subsubsection{Order $\eps^0$ and $\eps$}

Inserting the ansatz~\eqref{ansatz} into~\eqref{hj2} and expanding in terms of 
$\eps,$ the $O(\eps^0)$ and $O(\eps^1)$ equations are
\begin{equation*}
\begin{split}
O(1):& \quad \partial_t S_0 - \partial_t \tau \ \partial_\beta S_1=
\frac{1}{2} P_1^T P_1 + \Vhat + (a - \partial_\beta S_1 - \partial_\gamma
S_1) \what, 
\\
O(\eps):& \quad \partial_t S_1 - \partial_t \tau \ \partial_\beta S_2 =
\nabla_1 \Vhat \cdot (\partial_{P_1} S_1 - \partial_{P_1} 
\tau \partial_\beta S_2) 
\\
&\qquad + \nabla_2 \Vhat \cdot ((x+\partial_Y S_1) \cos \gamma + Y \sin
\gamma) 
\\
&\qquad + (a - \partial_\beta S_1
 - \partial_\gamma S_1) \nabla_1 \what \cdot 
(\partial_{P_1} S_1 - \partial_{P_1} \tau \ \partial_\beta S_2) 
 - (\partial_\beta S_2 + \partial_\gamma S_2) \what, 
\end{split}
\end{equation*}
where 
$$
\Vhat = V(q_1+\partial_{P_1} S_0 - \partial_{P_1} \tau \ 
\partial_\beta S_1,0) 
\quad \text{and} \quad
\what = \w(q_1 + \partial_{P_1} S_0 - \partial_{P_1} \tau \ 
\partial_\beta S_1).
$$  
As $\Vhat$ and $\what$ depend on $S_1,$ we must proceed carefully in 
closing the equations.  In fact, we will show $S_1 =0$ and close the
equations on $S_0$ by expanding
the unknowns $S_0$ and $S_1$ in powers of $t.$

\begin{lemma}\label{thm:split}
The $O(1)$ and $O(\eps)$ equations imply
\begin{align}
\partial_t S_0 &= \frac{1}{2} P_1^T P_1 + V(q_1+ \partial_{P_1} S_0, 0) 
+ a \w(q_1+\partial_{P_1} S_0),
\label{S0}
\\
S_1 &= 0. \label{S1} 
\end{align}
\end{lemma}

\begin{proof}
We expand the functions $S_0$, $S_1$ and $\tau$ in powers of $t$ : 
$$
S_0 = S_{0,0} + t S_{0,1} + \frac{t^2}{2} S_{0,2} + \cdots,
\ \
S_1 = S_{1,0} + t S_{1,1} + \frac{t^2}{2} S_{1,2} + \cdots,
\ \
\tau = \tau_0  + t \tau_1  + \frac{t^2}{2} \tau_2  + \cdots
$$
where from the initial conditions of~\eqref{tau} and~\eqref{ICs}, 
we have that $S_{0,0} = \tau_0 =0.$
The $O(\eps^0 t^0)$ equation is
\begin{equation*}
S_{0,1} - \tau_1 \partial_\beta S_{1,0} = \frac{1}{2} P_1^T P_1 +
V(q_1,0)
+ (a - \partial_\beta S_{1,0} - \partial_\gamma S_{1,0}) \w(q_1).
\end{equation*}
We wish to separate the $S_0$ and $S_1$ terms. 
We first observe that~\eqref{tau} gives $\tau_1 = \w(q_1)$ which
leaves 
\begin{equation*}
S_{0,1} = \frac{1}{2} P_1^T P_1 + V(q_1,0) 
            + (a - \partial_\gamma S_{1,0}) \w(q_1).
\end{equation*}
We integrate both sides of the equation with respect to $\gamma$ between
0 and $2\pi$ and note that, by the periodicity assumption in the ansatz, 
$\dps \dashint_0^{2\pi} \partial_\gamma S_{1,0} \
d\gamma = 0$.
Since $S_0$ is independent of $\gamma$, this leaves 
\begin{equation*}
S_{0,1} = \frac{1}{2} P_1^T P_1 + V(q_1,0) + a \w(q_1), 
\end{equation*}
which is~\eqref{S0} up to $O(t)$ error.
This implies $\partial_\gamma S_{1,0} = 0$ which, combined with the initial 
condition on $S_1,$ implies $S_{1,0} = 0,$ so that~\eqref{S1} holds up to
$O(t)$ error.
  
From here we proceed by induction.  Suppose that $\eqref{S0}$
and $\eqref{S1}$ hold up to
$O(t^k).$  We note from~\eqref{tau}
that 
\begin{equation*}
\partial_{P_1} \tau = O(t^2),
\end{equation*}
and thus $\partial_{P_1} \tau \ \partial_\beta S_1 = O(t^{k+2}).$
We then infer from the $O(1)$ equation that 
\begin{equation*}
\begin{split}
\partial_t S_0 - \frac{t^{k}}{k!} \tau_1 \partial_\beta S_{1,k} 
&= \frac{1}{2} P_1^T P_1 + V(q_1+\partial_{P_1} S_0,0) 
+ a \w(q_1+\partial_{P_1} S_0)
\\
&\qquad - \frac{t^{k}}{k!} 
(\partial_\gamma S_{1,k} + \partial_\beta
S_{1,k}) \w(q_1) + O(t^{k+1}).
\end{split}
\end{equation*}
Upon canceling out the two $\partial_\beta S_{1,k}$ terms using $\tau_1 = \w(q_1)$ 
and integrating with respect to $\gamma,$ we find
that~\eqref{S0} holds up to $O(t^{k+1})$ and
that $\partial_\gamma S_{1,k} =0.$

Because $S_1$ vanishes up to an error of $O(t^k),$ 
we infer from~\eqref{tau} that  
\begin{equation*}
\partial_{t} \tau - \big(\nabla_1 \Vhat + a \nabla_1 \what\big) 
\cdot \partial_{P_1} \tau
- \what = O(t^{k+2}).
\end{equation*}
We first use this in the $O(\eps)$ equation to cancel the
$\partial_\beta S_2$ terms. We then neglect all terms of $O(t^k)$
to find
\begin{equation*}
\frac{t^{k-1}}{(k-1)!} S_{1,k} =\nabla_2 V(q_1+\partial_{P_1} S_0,0) 
\cdot (x \cos\gamma + Y \sin\gamma)
- \w(q_1 + \partial_{P_1} S_0) \partial_\gamma S_{2} + O(t^{k}).
\end{equation*}
Integrating this equation with respect to $\gamma$ gives $S_{1,k} = 0.$  Thus,
we have shown that~\eqref{S1} holds up to $O(t^{k+1})$. Consequently, we
have shown by induction that~\eqref{S0} and~\eqref{S1} hold to all orders.
\end{proof}

In view of~\eqref{S1}, we have the following simplified
expressions for $\Vhat$ and $\what$:
$$
\Vhat = V(q_1+\partial_{P_1} S_0,0) 
\quad \text{and} \quad
\what = \w(q_1 + \partial_{P_1} S_0).
$$
In the proof to Lemma~\ref{thm:split} we have found 
$S_{0,1} = \frac{1}{2} P_1^T P_1 + V(q_1,0) + a \w(q_1).$
Using~\eqref{S0}, we likewise find 
$S_{0,2} = \nabla_1 V(q_1,0) \cdot P_1 +  a \nabla_1 \w(q_1) \cdot P_1.$  
For our generating function, we keep these two lowest order
terms. We combine them, neglecting only $O(t^3)$ terms, to give
\begin{equation}
\label{S0app}
\widetilde{S}_0 = t \left[\frac{1}{2} P_1^T P_1 + V\Big(q_1 + \frac{t}{2} P_1,0\Big) 
+ a \w\Big(q_1 + \frac{t}{2} P_1\Big) \right].
\end{equation}
We likewise keep the first two orders in the fast time $\tau.$
From~\eqref{tau}, we have $\tau_0 = 0,$ $\tau_1 = \w(q_1),$ and
$\tau_2 = \nabla \w(q_1) \cdot P_1.$
We combine these, neglecting only $O(t^3)$ terms, and take the approximation
\begin{equation}
\label{tauapp}
\widetilde{\tau} = t \w\Big(q_1+\frac{t}{2} P_1\Big). 
\end{equation}
Inserting $S_1 = 0$ into the $O(\eps)$ equation and using 
equation~\eqref{tau} leaves
\begin{equation}
\what \ \partial_\gamma S_2 =  \nabla_2 \Vhat \cdot (x \cos \gamma 
+ Y \sin \gamma).
\label{S2g}
\end{equation}

%%%%%%%%%%%%%%%%%%%%%%%%%%%%%%%%%%%%%%%%%%%%%%%%%%%%%%%%%%%%%%%%%%%%%%%%%
\subsubsection{Order $\eps^2$}
The $O(\eps^2)$ equation of~\eqref{hj2} is 
\begin{equation}
\begin{split}
\label{findS2}
\partial_t S_2 - \partial_t \tau \ \partial_\beta S_3 
&= \frac{1}{2} ( \nabla_{11} \Vhat + a  \nabla_{11} \what) :
                   (\partial_{P_1} \tau \otimes \partial_{P_1} \tau)
 (\partial_\beta S_2)^2 \\
&\qquad+ (\nabla_1 \Vhat +a \nabla_1 \what) \cdot 
    (\partial_{P_1} S_2 - \partial_{P_1} \tau \ \partial_\beta S_3) \\
&\qquad- \nabla_{12} \Vhat : (\partial_\beta S_2 \ \partial_{P_1} \tau 
                         \otimes \ell) 
 + \frac{1}{2} \nabla_{22} \Vhat : \ell \otimes \ell
 + \nabla_2 \Vhat \cdot \partial_Y S_2 \cos \gamma \\
&\qquad + (\partial_\gamma S_2 + \partial_\beta S_2) \ 
             \nabla_1 \what \cdot (\partial_{P_1} \tau \ \partial_\beta S_2)
 - (\partial_\gamma S_3 + \partial_\beta S_3) \ \what
\end{split}
\end{equation}
where $\ell = x \cos \gamma + Y \sin \gamma.$
Since our approximation $\widetilde{\tau}$ given in~\eqref{tauapp} 
satisfies~\eqref{tau} up to $O(t^2),$
we can cancel the terms involving
$\partial_\beta S_3$ up to $O(t^2),$ which leaves 
\begin{equation*}
\begin{split}
\partial_t S_2 
&= \frac{1}{2} ( \nabla_{11} \Vhat + a  \nabla_{11} \what) :
                   (\partial_{P_1} \widetilde{\tau} \otimes \partial_{P_1} \widetilde{\tau})
 (\partial_\beta S_2)^2 + (\nabla_1 \Vhat +a \nabla_1 \what) \cdot 
    \partial_{P_1} S_2  \\
&\qquad- \nabla_{12} \Vhat : (\partial_\beta S_2 \ \partial_{P_1} \widetilde{\tau} 
                         \otimes \ell) 
 + \frac{1}{2} \nabla_{22} \Vhat : \ell \otimes \ell
 + \nabla_2 \Vhat \cdot \partial_Y S_2 \cos \gamma \\
&\qquad + (\partial_\gamma S_2 + \partial_\beta S_2) \ 
             \nabla_1 \what \cdot (\partial_{P_1} \widetilde{\tau} \ \partial_\beta S_2)
 - \partial_\gamma S_3 \ \what + O(t^2).
\end{split}
\end{equation*}
As before, we integrate with respect to $\gamma$ between 0 and
$2\pi$ to remove the $\partial_\gamma S_3$ term, leaving only an equation in
$S_2.$  This time, however, $S_2$ is not itself independent of $\gamma.$  
From~\eqref{S2g} we can write 
\begin{equation*} 
\begin{split} 
S_2(t,\beta,q_1,x,\gamma,P_1,Y,a) &= \frac{\nabla_2 \Vhat}{\what}
\cdot (x \sin \gamma - Y \cos \gamma)
+ C(t,\beta,q_1,x,P_1,Y,a).
\end{split} 
\end{equation*} 
The initial condition $S_2(0,\gamma,q_1,x,\gamma,P_1,Y,a)=0$
implies that 
\begin{equation}
\label{C0}
C(0,\beta,q_1,x,P_1,Y,a) = - \frac{\nabla_2 V(q_1,0)}{\w(q_1)}
\cdot (x \sin \beta - Y \cos \beta).
\end{equation}
We integrate~\eqref{findS2} with respect to $\gamma$ between 0 and
$2\pi$, which gives
\begin{equation}
\begin{split}
\label{findC}
 \partial_t C  &= 
 \nabla_{22} \Vhat :  \frac{(x \otimes x + Y \otimes Y)}{4}
- \frac{(\nabla_2 \Vhat)^2}{2 \what} 
+(\nabla_1 \Vhat + a \nabla_1 \what) \cdot \partial_{P_1} C \\
&\qquad + \Big[ \nabla_1 \what \cdot \partial_{P_1} \widetilde{\tau} + 
               \frac{1}{2} (\nabla_{11} \Vhat 
+ a \nabla_{11} \what) : (\partial_{P_1} \widetilde{\tau} \otimes \partial_{P_1} \widetilde{\tau}) \Big] \
(\partial_\beta C)^2  + O(t^2)
\end{split}
\end{equation}
where we have used
\begin{equation*}
\partial_t C = \dashint_0^{2\pi} \partial_t S_2 \ d \gamma,
\qquad
\partial_\beta C = \dashint_0^{2\pi} \partial_\beta S_2 \ d \gamma,
\qquad
\partial_{P_1} C = \dashint_0^{2\pi} \partial_{P_1} S_2 \ d \gamma,
\end{equation*}
\begin{equation*}
\begin{split} 
\dashint_0^{2 \pi} \partial_Y S_2 \cos \gamma \ d \gamma
&= \dashint_0^{2 \pi} \left[ - \frac{\nabla_2 \Vhat}{\what} \cos^2 \gamma 
+ \partial_Y C \cos\gamma \right] \ d \gamma =  - \frac{\nabla_2 \Vhat}{2\what},
\end{split}
\end{equation*}
and
\begin{equation*}
\begin{split} 
\dashint_0^{2 \pi} 
   \frac{1}{2} \nabla_{22} \Vhat : \ell \otimes \ell \ d\gamma 
&= \frac{1}{2} \nabla_{22} \Vhat: \dashint_0^{2 \pi} 
     (x \otimes x \cos^2 \gamma + x \otimes Y \cos \gamma \sin \gamma 
        + Y \otimes Y \sin^2 \gamma) \  d\gamma \\
&= \frac{1}{4} \nabla_{22} \Vhat : (x \otimes x + Y \otimes Y).
\end{split}
\end{equation*} 
We expand 
$$
C(t,\beta,q_1,x,P_1,Y,a) = C_0(\beta,q_1,x,P_1,Y,a) 
 + t C_1(\beta,q_1,x,P_1,Y,a) + \frac{t^2}{2} C_2(\beta,q_1,x,P_1,Y,a) +
\cdots
$$ 
and insert this expansion into~\eqref{findC}. The $O(t)$ term satisfies 
\begin{equation*}
C_1 = 
 \frac{1}{4} \nabla_{22} V(q_1,0) : (x \otimes x + Y \otimes Y)
- \frac{(\nabla_2 V(q_1,0))^2}{2 \w(q_1)}.
\end{equation*} 
We keep the $O(t^0)$
term~\eqref{C0} along with the $O(t)$ term to approximate $S_2$ up to $O(t^2)$ error with 
\begin{equation}
\begin{split}
\widetilde{S}_2 &= 
 \frac{\nabla_2 V(q_1+t P_1,0)}{\w(q_1+t P_1)}
\cdot (x \sin \gamma - Y \cos \gamma)
- \frac{\nabla_2 V(q_1,0)}{\w(q_1)}
\cdot (x \sin \beta - Y \cos \beta) \\
&+ t \left[ \frac{1}{4} \nabla_{22} V(q_1,0) : (x \otimes x + Y \otimes Y)
- \frac{(\nabla_2 V(q_1,0))^2}{2 \w(q_1)} \right].
\end{split}
\label{S2app}
\end{equation} 

\subsection{Generating function and algorithm}
\label{sec:finite_diff}

Combining our previous approximations~\eqref{S0app} and~\eqref{S2app},
we approximate the solution to~\eqref{hj2} by 
\begin{equation}
\label{gen}
\begin{split}
\widetilde{S}_\eps &= \widetilde{S}_0 + \eps^2 \widetilde{S}_2 
\\
&=t \left[ \frac{1}{2} P_1^T P_1 + V\left(q_1 + \frac{t}{2} P_1, 0\right) 
       + a \w\left(q_1 + \frac{t}{2} P_1\right) \right]
\\
&+ \eps^2 \left[ \frac{\nabla_2 V(q_1+tP_1,0)}{\w(q_1+tP_1)} 
                           (x \sin \Sigma - Y \cos \Sigma) \right.\\
&
- \frac{\nabla_2 V(q_1,0)}{\w(q_1)} \left(x \sin \left(\Sigma - \frac{t}{\eps}
\w(q_1+\frac{t}{2}P_1) \right) 
                               - Y \cos \left(\Sigma - \frac{t}{\eps}
\w(q_1+\frac{t}{2}P_1) \right) \right)
\\
&\left.+ \frac{t}{4} \nabla_{22} V(q_1,0) : (x \otimes x+ Y \otimes Y) 
- t\frac{(\nabla_2 V)^2}{2 \w}
\right],
\end{split}
\end{equation}
where we note 
$$
S_\eps(t) = \widetilde{S}_\eps + O(t^3) + O(\eps^2 t^2) + O(\eps^3).
$$

The existence of derivatives of $V$ in the generating function is
problematic because it leads to higher-order gradient computations in
the resulting numerical scheme.  In many cases of interest, these
higher-order gradients are computationally expensive.  Furthermore,
computing them potentially involves extra effort during
implementation as most traditional schemes only involve first
derivatives of $V$.
To avoid this issue in our numerical
scheme, we replace all derivatives of $V$ found in the generating function with
finite differences.
We make this replacement before inserting the generating function
into~\eqref{symp_scheme} in order to
maintain the symplecticity of the scheme.  Note that there is no unique way to
form the finite differences.  

In doing this, we choose to exclude  $\dps - \frac{(\nabla_2 V)^2}{2 \w}$ from the
$O(\eps^2 t)$ term as it is much more computationally expensive to
implement as a finite difference than the other terms.  
Rather, we only retain a portion of the 
$O(\eps^2 t)$ term that turns out to be sufficient to reproduce 
the exchange of fast
actions (see the definition \eqref{Ij} and simulation in
Fig.~\ref{fig:traj} below).  Note that the only change in the order of
accuracy is caused by not retaining the portion of the $O(\eps^2 t)$
term mentioned above.  We choose the approximation
\begin{equation*}
\begin{split}
\SFD  &= t \left[ \frac{1}{2} P_1^T P_1 + V \left(q_1 +
\frac{t}{2} P_1, 0\right) + 
a \w\left( q_1+\frac{t}{2}P_1 \right) \right] 
\\
&+ \frac{\eps}{\w(q_1+tP_1)}
\bigg[ V(q_1+tP_1, \eps ( x \sin \Sigma - Y \cos \Sigma))
- V(q_1+tP_1,0) \bigg] 
\\
&+ \frac{\eps}{\w(q_1)} \left[
V(q_1,0) - V\left(q_1, \eps \left( x \sin (\Sigma - \frac{t}{\eps}
\w(q_1+\frac{t}{2}P_1)) - Y \cos (\Sigma - \frac{t}{\eps}
\w(q_1+\frac{t}{2}P_1))\right)\right)
\right]
\\
&+ \frac{t}{4} \Big[ 
V(q_1,\eps x) - 2 V(q_1,0) + V(q_1,-\eps x)
+ V(q_1,\eps Y) - 2 V(q_1,0) + V(q_1,-\eps Y)
 \Big].
\end{split} 
\end{equation*}
Substituting $S_{\eps,{\rm FD}}$ into~\eqref{symp_scheme}, we arrive at
the following expressions, where we use 
$\dps \sigma = \Sigma - \frac{h}{\eps} \w \left(q_1 + \frac{h}{2}
  P_1\right)$. We first solve for $(P_1,Y,\Sigma)$ in the implicit equations
\begin{equation}
\label{alg1}
\begin{split}
P_1 &= p_1 - h [\nabla_1 V(q_1 + \smfrac{h}{2} P_1, 0)  
    +  a \nabla_1 \w (q_1 + \smfrac{h}{2} P_1) ]\\
&\qquad - \eps h \nabla_1 \w(q_1+\smfrac{h}{2}P_1) \frac{\nabla_2 \Vd}{\w(q_1)} 
          [x \cos \sigma + Y \sin\sigma ] \\
&\qquad - \eps \frac{\nabla_1 \Ve - \nabla_1 \Vc}{\w(q_1+hP_1)} \\
&\qquad + \eps \frac{\nabla_1 \w(q_1+hP_1) (\Ve - \Vc)}{\w(q_1+hP_1)^2} \\
&\qquad - \eps \frac{\nabla_1 \Va - \nabla_1 \Vd}{\w(q_1)} \\
&\qquad + \eps \frac{\nabla_1 \w(q_1) (\Va - \Vd)}{\w(q_1)^2} \\
&\qquad - \frac{h}{4} \Big(\nabla_{1} V(q_1,\eps x) -2 \nabla_{1} V(q_1,0)
 + \nabla_{1} V(q_1,-\eps x) \\
&\qquad\qquad+ \nabla_{1} V(q_1,\eps Y) -2 \nabla_{1} V(q_1,0)
 + \nabla_{1} V(q_1,-\eps Y) \Big),
\\ 
Y &= y - \eps \frac{\nabla_{2} \Ve}{\w(q_1+hP_1)} \sin \Sigma 
       + \eps \frac{\nabla_{2} \Vd}{\w(q_1)} \sin \sigma \\
&\qquad       - \frac{h}{4} \Big(\nabla_{2} V(q_1,\eps x) - 
                          \nabla_2 V(q_1,-\eps x)\Big),  
\\
\Sigma &= \sigma + \frac{h}{\eps} \w(q_1+\smfrac{h}{2}P_1).
\end{split} 
\end{equation}
We next compute $(Q_1,X,A)$ using 
\begin{equation}
\label{alg2}
\begin{split}
Q_1 &= q_1 + h P_1 + \smfrac{h^2}{2} \Big[\nabla_1 V(q_1 + \smfrac{h}{2} P_1, 0) 
    + a \nabla_1 \w (q_1 + \smfrac{h}{2} P_1) \Big]\\
&\qquad + \eps \smfrac{h^2}{2}  \nabla_1 \w(q_1 + \smfrac{h}{2}P_1)
\frac{\nabla_2 \Vd}{\w(q_1)} [x \cos \sigma + Y \sin\sigma ] \\
&\qquad+ h \eps \frac{\nabla_{1} \Ve - \nabla_1 \Vc}{\w(q_1+hP_1)} \\
&\qquad- h \eps \frac{(\Ve - \Vc) \nabla_1 \w(q_1+hP_1)}{\w(q_1+hP_1)^2}, 
\\
X &= x - \eps \frac{\nabla_{2} \Ve}{\w(q_1+hP_1)} \cos \Sigma 
       + \eps \frac{\nabla_{2} \Vd}{\w(q_1)} \cos \sigma \\
&\qquad       + \frac{h}{4} \Big(\nabla_{2} V(q_1,\eps Y) - 
                          \nabla_2 V(q_1,-\eps Y)\Big),  
\\
A &= a + \eps \frac{\nabla_2 \Vd}{\w(q_1)} 
            [x \cos \sigma + Y \sin\sigma ] \\
&\qquad  - \eps \frac{\nabla_{2} \Ve}{\w(q_1+hP_1)} 
            [x \cos \Sigma + Y \sin \Sigma ], 
\end{split} 
\end{equation}
where we recall $\dps \sigma = \Sigma -
\frac{h}{\eps} \w \left(q_1 + \frac{h}{2} P_1\right)$. 
The equations~\eqref{alg1} are implicit in $Z = (P_1,Y,\Sigma)$ and 
can be written in the form 
\begin{equation}
\label{implicit}
Z = z + h F(Z) + \eps G(Z) + \frac{h}{\eps} K(Z) 
= z + {\cal A}(Z),
\end{equation}
with $z=(p_1,y,\sigma)$ and where $\dps K(Z)=(0,0,\w(q_1+\frac{h}{2}P_1))$ 
(the dependency of $F$, $G$ and $K$ on $(q_1,x,a)$ is not explicitly
written).  We solve~\eqref{implicit} for $Z$ with a simple fixed point
iteration.  In our simulations the gradient of ${\cal A}$ with respect to
$Z$ is small 
since the parameters $h$ and $\eps$ are small and since in practice we work
with $h^2/\eps < 0.4$  (see e.g. Fig.~\ref{fig:res}
below).  We thus expect the
fixed point algorithm to quickly converge.  In the test case considered below, this is indeed the case.
We terminate the iteration when successive iterates differ by a
relative factor of $10^{-10}$ or less, and the algorithm typically converges
in 3 to 8 iterations, depending on the stepsize $h$ and the parameter
$\eps$. 

\medskip

The equations~\eqref{cov}, \eqref{alg1} and~\eqref{alg2} yield
Algorithm~\ref{algo:alg} outlined below.

\bookbox{
\begin{algorithm}
Initialize: From the initial conditions
$(\check{q}_1(0),\check{q}_2(0),\check{p}_1(0),\check{p}_2(0))$, compute
the initial conditions
\begin{align*}
q_1(0) &= \check{q}_1(0),  &
p_1(0) &= \check{p}_1(0) 
    - \frac{\nabla_1 \w(\check{q}_1(0))}{2 \w(\check{q}_1(0))} 
                         \check{q}_2(0)^T \check{p}_2(0),  \\
x(0) &= \frac{\sqrt{\w(\check{q}_1(0))}}{\eps} \check{q}_2(0), &
y(0) &= \frac{1}{\sqrt{\w(\check{q}_1(0))}} \check{p}_2(0),  \\
\sigma(0) &= 0, &
a(0) &= \frac{1}{2 \w(\check{q}_1)} \left(
       \check{p}_2^T(0) \check{p}_2(0) +
       \frac{\w(\check{q}_1)^2}{\eps^2} \check{q}_2^T(0) \check{q}_2(0) 
      \right).
\end{align*}
Iterate: for $n \geq 0$,
\begin{enumerate}
\item
Set $(q_1,x,\sigma,p_1,y,a) = 
\left(q_1^n,x^n,\sigma^n,p_1^n,y^n,a^n\right)$.
\item Solve the implicit equations~\eqref{alg1} for $(P_1,Y,\Sigma)$.
\item Compute $(Q_1,X,A)$ using~\eqref{alg2}.
\item Set $\left(q_1^{n+1},x^{n+1},\sigma^{n+1},p_1^{n+1},y^{n+1},a^{n+1}\right)
= (Q_1,X,\Sigma,P_1,Y,A)$. 
\end{enumerate}
Post-process:  From $\left(q_1^N,x^N,\sigma^N,p_1^N,y^N,a^N\right)$,
return to the original variables: compute first $(q_2^N,p_2^N)$ using
$$
q_2^N = \sqrt{\eps} \left[ x^N \cos(\sigma^N) + y^N \sin(\sigma^N) \right],
\quad 
p_2^N = \sqrt{\eps} \left[ -x^N \sin(\sigma^N) + y^N \cos(\sigma^N) \right],
$$
and compute next $(\check{q}_1^N,\check{q}_2^N,\check{p}_1^N,\check{p}_2^N)$ using
$$
\check{q}_1^N = q_1^N, 
\qquad
\check{p}_1^N = p_1^N + \frac{\nabla_1 \w(q_1^N)}{2 \w(q_1^N)} \ 
q_2^N \cdot p_2^N, 
\qquad 
\check{q}_2^N = \frac{\sqrt{\eps}}{\sqrt{\w(q_1^N)}} q_2^N,
\qquad
\check{p}_2^N = \frac{\sqrt{\w(q_1^N)}}{\sqrt{\eps}} p_2^N.
$$

\label{algo:alg}
\end{algorithm}
}

As a variant, we will also consider a ``no-loop'' version of Algorithm~\ref{algo:alg}.
Instead of solving~\eqref{implicit} for $Z$, we perform the
following update:
\begin{equation}
\label{noloop}
\begin{split}
z^\star &= z + h F(z) + \frac{h}{\eps} K(z),
\\
Z^\star &= z + h F(z^\star) + \eps G(z^\star) + \frac{h}{\eps} K(z^\star),
\end{split}
\end{equation}
and approximate the solution $Z=(P_1,Y,\Sigma)$ to~\eqref{implicit}
(namely, \eqref{alg1}) by
$Z^\star = (P_1^\star,Y^\star,\Sigma^\star)$. Using $Z^\star$, we next
determine $(Q_1,X,A)$ in an explicit fashion, using~\eqref{alg2}. This yields an explicit
scheme (we first determine $z^\star$, next $Z^\star$ and finally
$(Q_1,X,A)$), with a lower computational cost than~\eqref{alg1}-\eqref{alg2}. However,
this scheme is not symplectic since we do not 
fully solve~\eqref{symp_scheme}. In Fig.~\ref{fig:eff} below, we will observe
that the no-loop version has comparable error behavior with a reduced
computational cost, though it seems difficult to justify why the
energy and invariants are well preserved in this version. 

%%%%%%%%%%%%%%%%%%%%%%%%%%%%%%%%%%%%%%%%%%%%%%%%%%%%%%%%%%%%%%%%%%%%%%%%%
\subsection{Numerical results} 

In this section, we provide numerical tests of the behavior of our
algorithm.  As mentioned above, we are not interested in 
computing exact trajectories for the fast variables.  Rather, we
are interested in how well the algorithm preserves invariants of
the system as well as its prediction of quantities derived from
the fast variables.   
This system undergoes a slow exchange between the actions~\eqref{Ij}
of the fast degrees of freedom, and we test how well the algorithm
captures this exchange.  We also test the robustness of the
algorithm by testing for the appearance of resonances.

We compare our algorithm to a well-known integrator for highly
oscillatory dynamical systems, the Mollify algorithm~\cite{mollify},
which is a modification of the previously proposed Impulse
(also known as Verlet-I/r-RESPA)
algorithm~\cite{impulse1, impulse2}.  Both Impulse and Mollify
follow a kick/oscillate/kick pattern and incorporate the
slow forces $\nabla \Vt$ 
only at the ``kick'' steps, which are separated by a macro
time step that is large with respect to the shortest timescale in the
solution.  This time step is typically larger than $\eps$ and is thus
larger than the stable regime for Verlet.  For the
``oscillate'' step, these methods integrate the fast forces using a
stepsize that is small with respect to $\eps.$   The Mollify algorithm
differs from the Impulse algorithm in how it incorporates the forces at
the ``kick'' steps in order to improve the stability of the algorithm
in the face of resonances, an issue we discuss later.
These algorithms are
designed to minimize the number of evaluations of the slow force,
with the assumption that the ``oscillate'' step is cheaper, or
comparable in cost, to a single evaluation of the slow force.  In
our tests, we have used the Verlet scheme for the ``oscillate'' step
within Mollify, with inner stepsize equal to $\eps/100.$

%%%%%%%%%%%%%%%%%%%%%%%%%%%%%%%%%%%%%%%%%%%%%%%%%%%%%%%%%%%%%%%%%%%%%%%%%
\subsubsection{Modified FPU}

\begin{figure}[htbp]
\begin{center}
\includegraphics[width=8cm]{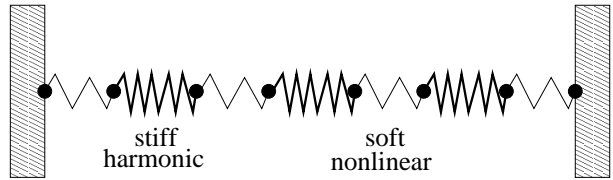}
\end{center}
\caption{\label{fig:fpu}Fermi-Pasta-Ulam spring chain.}
\end{figure}

The Fermi-Pasta-Ulam (FPU) chain is a commonly-used test case for highly
oscillatory integrators~\cite[Sec. XIII.2.1]{HLW}.  The chain is a
collection of alternating stiff, harmonic springs and soft, nonlinear
springs (see Fig.~\eqref{fig:fpu}).  
After a change of coordinates, the potential can be written
as in~\eqref{ham}, with $\w$ constant.  We
choose $f=s=3,$ corresponding to 3 stiff springs and 4 soft springs
(there is one less degree of freedom since the total chain length is
prescribed).  For positions $\check{q}_1 = (\check{q}_{1,1}, \check{q}_{1,2}, \check{q}_{1,3}) \in \RR^3$ and 
$\check{q}_2 = (\check{q}_{2,1}, \check{q}_{2,2}, \check{q}_{2,3}) \in \RR^3,$ the potential $\Vt$
is given by
\begin{equation*}
\Vt(\check{q}_1, \check{q}_2) =
\frac{1}{4} \left( (\check{q}_{1,1} - \check{q}_{2,1})^4 + 
\sum_{i=1}^2 
(\check{q}_{1,i+1} - \check{q}_{2,i+1} - \check{q}_{1,i} - \check{q}_{2,i})^4
+ (\check{q}_{1,3} + \check{q}_{2,3})^4 \right).
\end{equation*}
We modify the FPU problem to have {\em non-constant fast frequencies}.  
We consider a Hamiltonian of the form~\eqref{ham}, choose 
\begin{equation*}
\Omega(\check{q}_1) = \sqrt{1+\check{q}_{1,1}^2},
\end{equation*}
and use $\Vt$ above for the slow potential. 

The behavior of the modified FPU is qualitatively similar to the original.  
The fast and slow variables have timescales of $O(\eps)$ and $O(1),$
respectively. In addition there is a slow
exchange of the fast actions 
\begin{equation}
\label{Ij}
I_j = \frac{1}{2 \w(\check{q}_1)}  \left[\check{p}_{2,j}^2
+ \frac{\w(\check{q}_1)^2}{\eps^2} \check{q}_{2,j}^2 \right],
\qquad j=1,2,3,
\end{equation}  
over long time periods of $O(\eps^{-1}).$
The quantity
\begin{equation}
\label{eq:I_sum}
I = I_1 + I_2 + I_3
\end{equation}
is an adiabatic invariant (see~\cite{bornemann,schuette}
and~\cite[Theorem XIII.6.3]{HLW}):  
it is an $O(1)$ quantity that is
almost preserved over long times, with oscillations only 
of magnitude $O(\eps)$. 

For all numerical experiments, the initial conditions are 
\begin{equation*}
\check{q}_{1} = (1,0,0)^T,
\qquad
\check{p}_{1} = (1,0,0)^T,
\qquad
\check{q}_2 = (\eps,0,0)^T,
\qquad
\check{p}_2 = (1,0,0)^T.
\end{equation*}
For these values, we have
$H_\eps(0) = 2.5 + 3 \eps^2 + 0.5 \eps^4 \approx 2.5$
and $\dps I(0) = \frac{1}{2} \left(
\frac{1}{\sqrt{2}} + \sqrt{2} \right) \approx 1.06$.

%%%%%%%%%%%%%%%%%%%%%%%%%%%%%%%%%%%%%%%%%%%%%%%%%%%%%%%%%%%%%%%%%%%%%%%%%
\subsubsection{Preservation of invariants and exchange of actions}

\begin{figure}[htbp]
\centerline{
\input{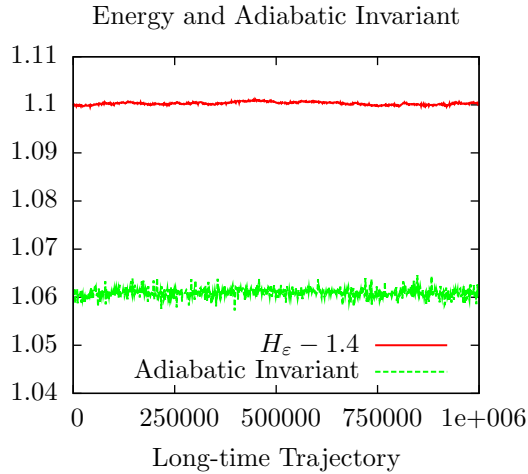} 
}
\caption{\label{fig:conserve} Energy (we plot $H_\eps - 1.4$ to allow
  better scaling) and adiabatic invariant computed with Algorithm~\ref{algo:alg}
over the long time interval $[0,10^6]$ for
$\eps=10^{-3}.$  We use stepsize $h=0.005.$ }
\end{figure}

We first monitor how the energy~\eqref{ham_encore} and the adiabatic
invariant~\eqref{eq:I_sum} are preserved along the numerical trajectory. 
Fig.~\ref{fig:conserve} shows the lack of drift in both quantities
over long times, for $h=0.005,$ $\eps=10^{-3}$ and
the above choice of initial conditions.  Note that energy preservation
for symplectic schemes is typically proven in the limit $h 
\rightarrow 0,$ for fixed $\eps.$  In our case, $h/\eps = 5,$ so that energy
preservation is not guaranteed by the theory.     

\medskip

\begin{figure}[tbp]
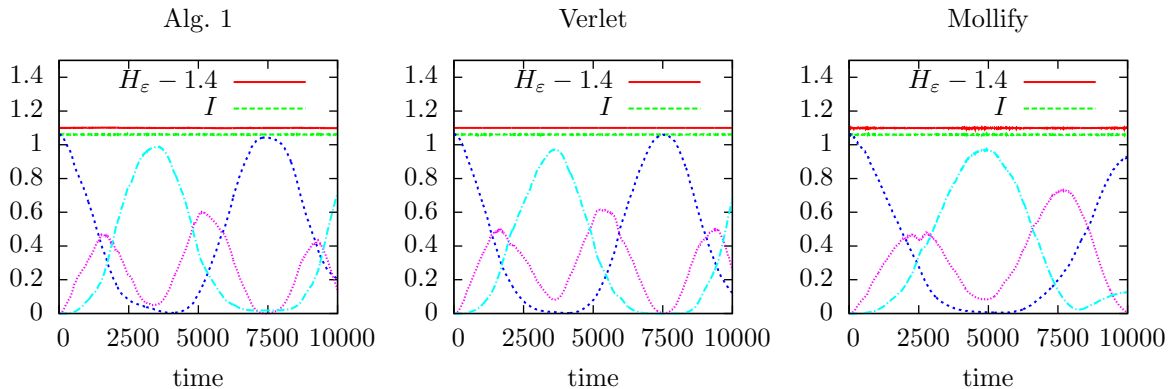

\centerline{
\begin{tabular}{ccc}
\input{figs/traj_M14_alg.tex} &
\hspace{-10mm}
\input{figs/traj_SV.tex} &
\hspace{-10mm}
\input{figs/traj_Moll.tex} 
\end{tabular}
}
\caption{\label{fig:traj}A single trajectory is simulated using three
algorithms, Algorithm~\ref{algo:alg}, Verlet, and Mollify.  The
energy, adiabatic invariant $I=I_1+I_2+I_3,$ and individual actions
$I_j$ are plotted
to time $10^4,$ with parameter $\eps=10^{-3}.$  
Verlet is run with a very small time step $h=10^{-5}$ and is
considered as the reference solution.  Mollify and
Algorithm~\ref{algo:alg} are run with
stepsizes $h=0.0008$ and $h=0.02,$ respectively, which are chosen so
that both algorithms involve roughly the same number of calls to the slow
forces.  This plot shows that Algorithm~\ref{algo:alg} better
captures the slow exchange of actions for a comparable computational
effort.}
\end{figure}

In Fig.~\ref{fig:traj}, we examine the slow exchange among the actions
$I_j$ defined in~\eqref{Ij}. We observe that Algorithm~\ref{algo:alg}
accurately reproduces exchange of the actions, in contrast to Mollify.
Note that we have used a smaller stepsize $h$ for Mollify than for
Algorithm~\ref{algo:alg} in order to balance one possible measure of
computational cost -- the number of evaluations of the slow force.  We
discuss the issue of computational efficiency in
Section~\ref{sec:cost}. 

%%%%%%%%%%%%%%%%%%%%%%%%%%%%%%%%%%%%%%%%%%%%%%%%%%%%%%%%%%%%%%%%%%%%%%%%%
\subsubsection{Resonance and computational efficiency}
\label{sec:cost}

\begin{figure}[tbp]
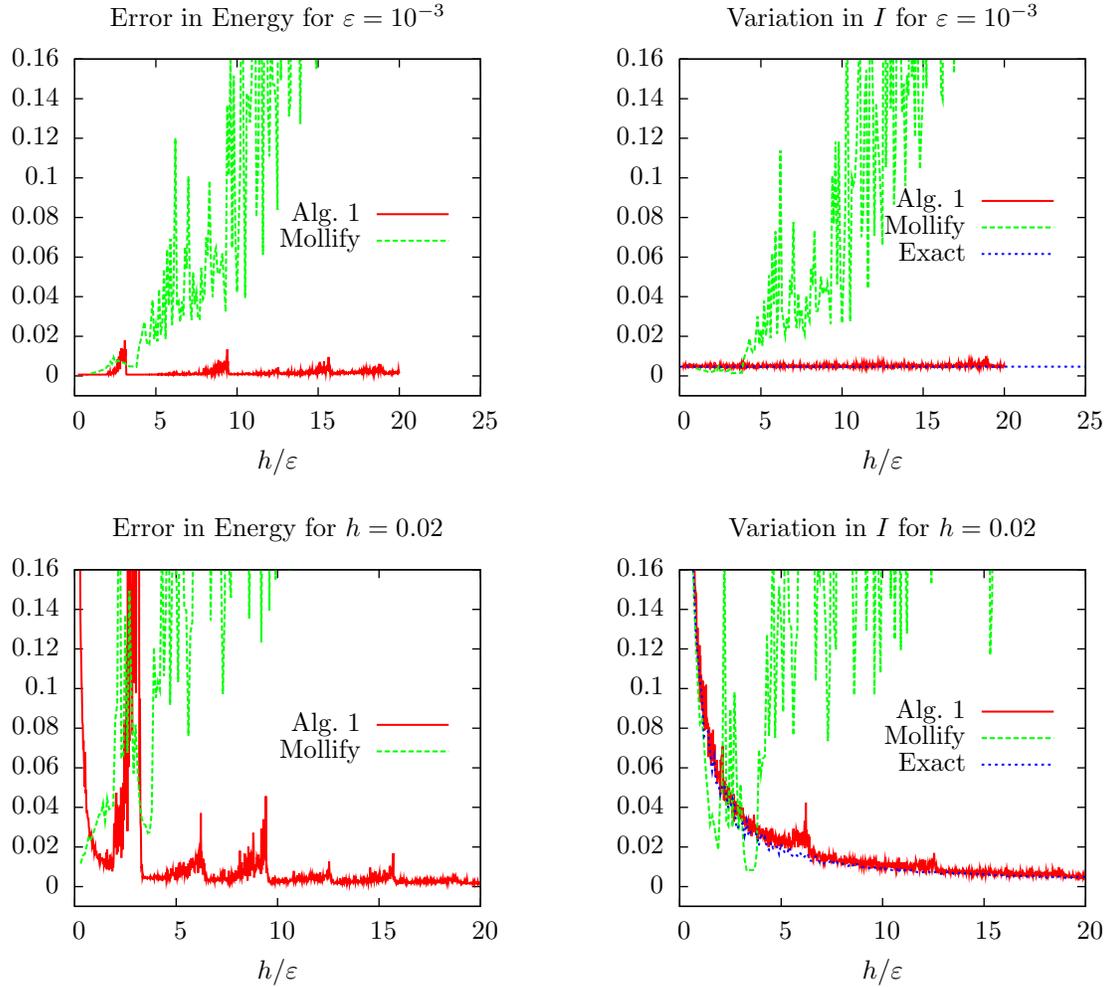

\centerline{
\begin{tabular}{cc}
\input{figs/hres_alg.tex} &
\input{figs/hres_ad_alg.tex} \\
\input{figs/eres_alg.tex} &
\input{figs/eres_ad_alg.tex}
\end{tabular}
}
\caption{\label{fig:res} Comparison of resonances for Algorithm~\ref{algo:alg} and
Mollify. In the upper row, a series of trajectories are
simulated for different stepsizes $h$, with $\eps=10^{-3}$.  In the
lower row, different values of $\eps$ are considered, and the dynamics is
integrated using the stepsize $h=0.02$. In both cases, the
trajectory is simulated to time $T=10^4.$ In the left (resp. right)
column, the maximum error in the energy (resp. the maximum variation in
$I$) over the trajectory, as defined in~\eqref{HI_err}, is plotted.  For the
variation of $I,$ a reference value is calculated using the Verlet
algorithm with a very small time step. Here, Mollify exhibits many more
resonances than it does in the case of constant $\w,$ as shown
in~\cite[Chap. XIII]{HLW} and~\cite[Figs. 4-10]{dcds}.} 
\end{figure}

One common difficulty encountered among integrators for highly
oscillatory systems is the appearance of resonances which destroy the
long-time preservation of energy and other invariants.  In linear
cases, these typically occur when the slow time step $h$ is an integer
multiple of half the period of the fast motion. In Fig.~\ref{fig:res},
we explore the resonance behavior by simulating many trajectories,
with various choices of the ratio $h/\eps.$  We do this in two
different ways: by holding $\eps$ fixed while varying $h$ as well as
by holding $h$ fixed while varying $\eps.$

We display the maximum error in energy and maximum variance in $I$ to time
$T=10^4,$
\begin{equation}
\label{HI_err}
{\rm err} = \max_{t \in [0,10^4]} |H_\eps(t) - H_\eps(0)|, \qquad
{\rm var} = \max_{t \in [0,10^4]} |I(t) - I(0)|.
\end{equation} 
Recall that $H_\eps(0) \approx 2.5$ and $I(0) \approx 1.06$.
Since the exact trajectory preserves energy, we desire that the method
has small error in $H_\eps.$  On the other hand, there is variation in
$I(t)$ even for the exact trajectory and correctly predicting the variation
is also of interest.  Therefore, for each $\eps$, we have also computed
a reference variation in $I$. Results are shown in Fig.~\ref{fig:res}. 

\medskip

For fixed $\eps=10^{-3}$ (the upper row of Fig.~\ref{fig:res}),
Algorithm~\ref{algo:alg} performs quite well in the whole range of $h,$ 
with only a few spikes in the error.  Mollify only performs well for small
$h,$ with large, generalized resonant regions. 
For fixed $h=0.02$ (the lower row of Fig.~\ref{fig:res}),
Algorithm~\ref{algo:alg} performs better for smaller $\eps,$ and the
error in the energy blows up for increasing $\eps.$  The prediction of
the variation in $I$ is even better, matching the reference result.
Note that the fact that Algorithm~\ref{algo:alg} performs 
better for smaller $\eps$ is consistent with the fact that it is derived
using homogenization techniques which are valid in the limit
$\eps\rightarrow 0.$ For large $\eps$, the 
terms neglected in the generating function are no longer small.   

\medskip

We next consider the maximum error in energy
versus the computational cost.  For the sake of our comparison here, we use
the number of evaluations of the slow force $-\nabla \Vt$ as a proxy for
overall computational cost.  We thus assume that the slow forces are
much more expensive to compute than the fast forces.
For example, in a simulation of molecular chains, if the fast forces 
represent the bonds between adjacent particles in the chain and the slow forces
include all other intramolecular as well as long-range intermolecular
forces, then the slow force is much more computationally involved. 
For the Mollify algorithm, this will mean assuming that the
computational expense of the `oscillate' step is negligible 
in comparison to evaluating the slow forces, which may, depending on
the application, underestimate the computational cost.  
The `oscillate' step involves propagating the position and momentum according
to the fast forces using a small time step and, in addition, propagating
the matrix $\left[ \begin{array}{cc} Q_q & Q_p \\ P_q &
P_p\end{array} \right]$
of partial derivatives of the position and momentum 
with respect to the initial condition~\cite{mollify,mollify2}.
This matrix of size $(s+f) \times (s+f)$ is propagated according the 
Hessian of the fast potential in the so-called variational equation.
In practice, simulating the variational equation is much more
expensive than
simulating a Hamiltonian dynamics on $q \in \RR^{s+f}$ according to the
fast forces.  In the sequel, the cost of the
`oscillate' step in Mollify is entirely neglected. 

\begin{figure}[tbp]
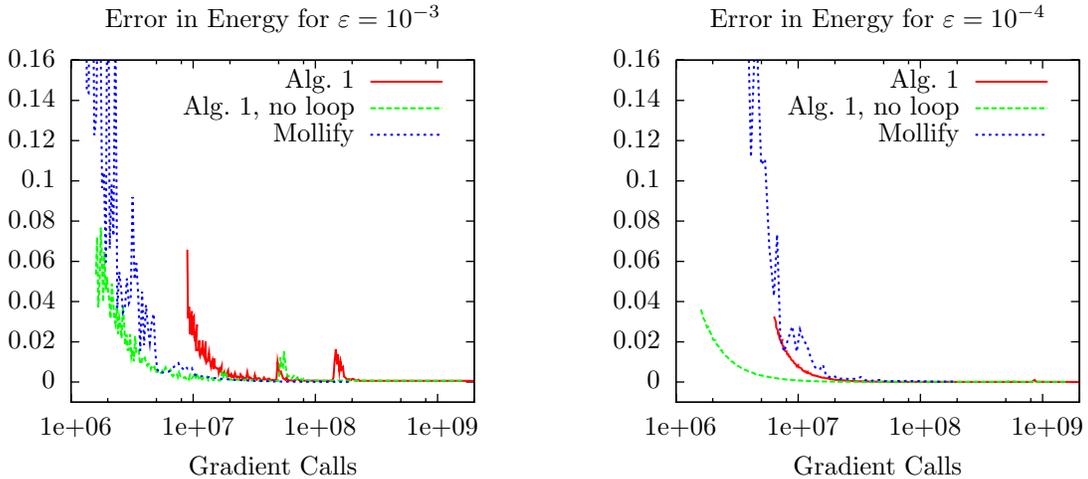

\centerline{
\begin{tabular}{cc}
\input{figs/eff3_alg.tex} &
\input{figs/eff4_alg.tex} %&
\end{tabular}
}
\caption{\label{fig:eff}Efficiency plots comparing
Algorithm~\ref{algo:alg}, a no-loop variant, and Mollify.  For each
graph, a series of trajectories are simulated for various stepsizes
$h$ with a single set of initial conditions and a single choice of
parameter $\eps$ ($\eps=10^{-3}$ for the left plot and $\eps=10^{-4}$
for the right plot).  The trajectory is simulated to time $T=10^4,$
and the maximum variation in the energy over the trajectory is
plotted.  
} 
\end{figure}

Figure~\ref{fig:eff} displays the maximum error in energy to time
$T=10^4$ versus the computational cost. 
Two plots are shown: we consider both cases $\eps = 10^{-3}$ and
$\eps = 10^{-4}$, and integrate the dynamics with different stepsizes
$h.$  As the stepsize decreases, the computational cost increases.  
For all algorithms, the error is generally decreasing as a function of 
computational cost, up to the presence of resonances.
In addition to Algorithm~\ref{algo:alg} and Mollify, we also consider
a variant of Algorithm~\ref{algo:alg}, that does not fully loop in
order to solve the implicit equations~\eqref{implicit} (see~\eqref{noloop}).

Although Algorithm~\ref{algo:alg} is much more stable with respect to
resonances (as shown in Fig.~\ref{fig:res}), the trade-off is that 
the computational cost per time-step is increased.
Thus, for $\eps=10^{-3}$, Mollify is slightly cheaper than
Algorithm~\ref{algo:alg} in terms of cost required to achieve a
certain accuracy. However, the
computational cost of the no-loop variant is much smaller than that of
Algorithm~\ref{algo:alg} and is also smaller than that of Mollify. 
In addition, the no-loop variant does not exhibit significantly more resonances
than does Algorithm~\ref{algo:alg}. 
For $\eps=10^{-4},$ Algorithm~\ref{algo:alg}
enjoys a slight edge in comparison to Mollify, and the no-loop variant
offers a significant advantage in computational cost.

We also observe from the comparison of the cases $\eps=10^{-3}$ and
$\eps=10^{-4}$ that, when $\eps$ decreases and $h$ is kept fixed, the
computational cost of Mollify remains the same, while its accuracy
decreases. On the other hand, the cost of Algorithm~\ref{algo:alg}
decreases somewhat since the number of iterations to solve the
implicit equations~\eqref{implicit} decreases, and the accuracy
increases because the terms neglected in the expansion become smaller.
In addition, resonances seem to disappear.

%%%%%%%%%%%%%%%%%%%%%%%%%%%%%%%%%%%%%%%%%%%%%%%%%%%%%%%%%%%%%%%%%%%%%%%%%
\section{The case of a matrix-valued, constant frequency}
\label{sec:multi}

In this section, we consider the case where the fast frequency is a
constant matrix. We consider the Hamiltonian
(see~\eqref{ham_matrix}) 
\begin{equation}
\label{ham_matrix2}
H_\varepsilon(q_1,q_2,p_1,p_2) = \frac{p_1^T p_1}{2} + \frac{p_2^T p_2}{2} +
V(q_1,q_2) + \frac{q_2^T \Omega^2 q_2}{2\varepsilon^2},
\end{equation}
and assume that $\Omega$ is a {\em
diagonal, positive definite} constant matrix of size $f \times f$
(if $\Omega$ is a constant positive definite symmetric matrix, then, upon
diagonalizing $\Omega$, we recover the case considered here). 
Our aim is to show that, upon a slight modification, we can also 
apply our strategy, thereby extending our previous study in~\cite{dcds}
which was restricted to the scalar case $\Omega = \omega \text{Id}$. In
the sequel, we have to distinguish two cases, whether or not the
eigenvalues of $\Omega$ are resonant. 

In both cases, we follow the exact same strategy as that we used
in~\cite{dcds} and in Section~\ref{sec:approach}: we first precondition
the fast motion using a change of variables and next apply our
two-scale ansatz to the Hamilton-Jacobi equation associated to the new
Hamiltonian. The difference between the non-resonant case and the
resonant case lies in the ansatz we make. In the non-resonant case, we
introduce a fast time for each frequency in the generating function
(see~\eqref{ansatz2_multi}). In contrast, in the resonant case, we
introduce a unique fast time for each group of frequencies that
are resonant one with each other (see~\eqref{ansatz2_reso}). In both
cases, the identification 
process follows the same lines. The derivation of the scheme is
presented in Sections~\ref{sec:non_reso} and~\ref{sec:reso}, in the
non-resonant case and in the resonant case, respectively. Numerical
results illustrating the non-resonant case are reported in
Section~\ref{sec:non_reso_num}. We conclude this section by considering
a test-case with three distinct frequencies, two of which are resonant
with each other (see Section~\ref{sec:multi_gene_num}).  

\subsection{The non-resonant case: derivation of the scheme}
\label{sec:non_reso}

Without loss of generality, we assume that the matrix $\Omega$
in~\eqref{ham_matrix2} reads
\begin{equation}
\label{eq:def_Omega_diago}
\Omega = \left(
\begin{array}{cccc}
\omega_1 \Id_{m_1} & 0 & \ldots & 0
\\
0 & \omega_2 \Id_{m_2} & \ldots & 0
\\
\ldots & \ldots & \ldots & \ldots
\\
0 & 0 & \ldots & \omega_d \Id_{m_d}
\end{array}
\right),
\end{equation}
with $0 < \omega_1 < \ldots < \omega_d$. The multiplicity of $\omega_i$
is $m_i$, with $\sum_{i=1}^d m_i = f$. 
In view of the block diagonal decomposition of $\Omega$, we
decompose $q_2 \in \RR^f$ into $d$ vectors $q_{2,i}  \in \RR^{m_i}$,
such that $q_2=(q_{2,1},\ldots,q_{2,d})$ and 
$$
q_2^T \Omega^2 q_2 = \sum_{i=1}^d \omega_i^2 \ q_{2,i}^T \, q_{2,i}.
$$
Likewise, we write $p_2=(p_{2,1},\ldots,p_{2,d})$.
We assume in the following that the $\omega_i$ are non-resonant, in the sense
that, for any $k \in \ZZ^d$,
\begin{equation}
\label{eq:non_reso}
\sum_{j=1}^d \omega_j k_j = 0 \Longrightarrow k = 0.
\end{equation}

\medskip

As pointed out above, we follow here the exact same strategy as that we
used in~\cite{dcds}, or in Section~\ref{sec:approach}: we first
precondition the fast motion, and next consider the Hamilton-Jacobi form
of the equations.
We proceed first with the time-dependent change of variables 
$(q_{2,i},p_{2,i}) \in \RR^{2 m_i} \mapsto
(x_{2,i},y_{2,i}) = \chi_i(t,q_{2,i},p_{2,i}) \in \RR^{2 m_i} $ defined by
$$
q_{2,i} = \cos \left(\frac{\omega_i t}{\varepsilon}\right) x_{2,i} +
\frac{\varepsilon}{\omega_i} \sin \left(\frac{\omega_i
    t}{\varepsilon}\right) y_{2,i} ,
\quad
p_{2,i} = - \frac{\omega_i}{\varepsilon} \sin \left(\frac{\omega_i
    t}{\varepsilon}\right) x_{2,i} + \cos \left(\frac{\omega_i
    t}{\varepsilon}\right) y_{2,i} ,
$$
which reads, in a more compact form,
\begin{equation}
\label{eq:cov_multi}
q_2 = \dps
\cos \left(\frac{\Omega t}{\varepsilon}\right) x_2 +
\varepsilon \Omega^{-1} \sin \left(\frac{\Omega
    t}{\varepsilon}\right) y_2 ,
\quad
p_2 = \dps
- \frac{\Omega}{\varepsilon} \sin \left(\frac{\Omega
    t}{\varepsilon}\right) x_2 + \cos \left(\frac{\Omega
    t}{\varepsilon}\right) y_2,
\end{equation}
with $x_2 = (x_{2,1},\ldots,x_{2,d})$ and $y_2 =
(y_{2,1},\ldots,y_{2,d})$. 
The dynamics on $(x_2,y_2)$ reads
\begin{eqnarray*}
\dot{x}_2 &=&
\varepsilon \Omega^{-1} \sin \left(\frac{\Omega
    t}{\varepsilon}\right) \partial_2 V \left[q_1(t),
\cos \left(\frac{\Omega t}{\varepsilon}\right) x_2(t) +
\varepsilon \Omega^{-1}   
\sin \left(\frac{\Omega t}{\varepsilon}\right) y_2(t) \right],
\\
\dot{y}_2
&=& - \cos \left(\frac{\Omega t}{\varepsilon}\right)
\partial_2 V \left[q_1(t), 
\cos \left(\frac{\Omega t}{\varepsilon}\right) x_2(t) +
\varepsilon \Omega^{-1}  
\sin \left(\frac{\Omega t}{\varepsilon}\right) y_2(t) \right],
\end{eqnarray*}
and 
the dynamics on $(q_1,x_2,p_1,y_2)$ is a Hamiltonian
dynamics with the time-dependent Hamiltonian
\begin{equation}
\label{eq:ham_non_reso}
H_\varepsilon^{\rm non-reso}(t,q_1,x_2,p_1,y_2) = 
\frac{p_1^T p_1}{2} + W_\eps^{\rm non-reso} \left(\frac{\omega_1
  t}{\eps},\ldots,\frac{\omega_d
  t}{\eps},q_1,x_2,y_2 \right), 
\end{equation}
where
$$
W_\eps^{\rm non-reso}(\tau_1,\ldots,\tau_d,q_1,x_2,y_2)
= V \left[q_1,
(\cos \tau_1) x_{2,1} 
+ \frac{\varepsilon}{\omega_1} (\sin \tau_1) y_{2,1}, \ldots,
(\cos \tau_d) x_{2,d} 
+ \frac{\varepsilon}{\omega_d} (\sin \tau_d) y_{2,d}
\right].
$$
We take the Hamiltonian~\eqref{eq:ham_non_reso} as a starting point for
our manipulations. 
In the sequel, we proceed with the construction of an approximation
$\widetilde{S_\eps}(h,q_1,x_2,P_1,Y_2)$ of the solution
$\overline{S_\eps}(h,q_1,x_2,P_1,Y_2)$ to
the Hamilton-Jacobi equation associated to~\eqref{eq:ham_non_reso},
for small times $h$. Observing that
the variable $x_2$ is of order $\eps$, we first perform a change of variables
and of unknown function:
$$
\alp_2 = \frac{\Omega}{\eps} x_2 \quad \text{and} \quad 
S_\eps(t,q_1,\alp_2,P_1,Y_2) = 
\overline{S}_\eps\left(t,q_1,\eps \Omega^{-1} \alp_2,P_1,Y_2\right),
$$
so that $S_\eps$ satisfies
\begin{equation}
\label{hj3_multi}
\left\{
\begin{array}{l}
\dps
\partial_t S_\eps = \frac{P_1^T P_1}{2} +
W_\eps^{\rm non-reso} \left(\frac{\omega_1 t}{\eps},\ldots,
\frac{\omega_d t}{\eps},q_1 + \partial_{P_1} S_\eps,
\eps \Omega^{-1} \alp_2 + \partial_{Y_2} S_\eps, Y_2 \right),
\\ 
S_\eps(0,q_1,\alp_2,P_1,Y_2) = 0.
\end{array}
\right.
\end{equation}
We make the ansatz
\begin{eqnarray}
\label{ansatz2_multi}
S_\eps(t,q_1,\alp_2,P_1,Y_2) &=& 
S_0 \left( t,\tau_1,\ldots,\tau_d,q_1,\alp_2,P_1,Y_2 \right) 
+
\eps S_1 \left( t,\tau_1,\ldots,\tau_d,q_1,\alp_2,P_1,Y_2 \right) 
\\
\nonumber
&+& \text{higher order terms in $\eps^k$, $k \geq 2$},
\end{eqnarray}
where the fast times $\tau_i$ are defined by
\begin{equation}
\label{eq:tau_multi}
\tau_i  = \frac{t \omega_i}{\eps}, \quad 1 \leq i \leq d,
\end{equation}
and where the functions $(S_k)_{k \geq 0}$ are assumed to be $2 \pi$ 
periodic with respect to each $\tau_i$. We also assume that the functions
$S_k$ are of
class $C^{1+\alpha}$ with respect to each $\tau_i$, with $\alpha \in \NN$,
$\alpha > d$.  

\begin{remark}
In the case where $V$ does not depend on $q_2$, the solution
to~\eqref{hj3_multi} can be analytically identified and is indeed of
the form~\eqref{ansatz2_multi}. 
\end{remark}

We now insert (\ref{ansatz2_multi}) in (\ref{hj3_multi}), identify the
first $d$ variables of $W_\eps^{\rm non-reso}$ with the fast times
$\left\{ \tau_i \right\}_{1 \leq i \leq d}$, and expand in powers of $\eps$. 
Based on~\eqref{symp_scheme}, we have that 
$X_2 = x_2 + \partial_{Y_2} \overline{S}_\eps$. 
The fast position $X_2$ is of order $\eps$, so $S_0$ does not depend on
$Y_2$. The equation of order $\eps^{-1}$ in the expansion
of~\eqref{hj3_multi} then becomes 
$$
\sum_{i=1}^d \omega_i \ \partial_{\tau_i} S_0 = 0.
$$
Using Lemma~\ref{lem:periodic} below, we deduce that 
$$
\forall 1 \leq i \leq d, \quad \partial_{\tau_i} S_0 = 0.
$$
\begin{lemma}
\label{lem:periodic}
Let $f(\tau_1,\ldots,\tau_d)$ be a function that is $2 \pi$ 
periodic with respect to each $\tau_j$, and of class $C^{1+\alpha}$
with respect to each of its argument, with $\alpha \in \NN$, $\alpha >
d$. Assume that  
\begin{equation}
\label{eq:assump}
\sum_{j=1}^d \omega_j \, \partial_{\tau_j} f = c,
\end{equation}
for a constant $c$, and a $d$-tuple $\left\{ \omega_j \right\}_{1 \leq j \leq d}$ satisfying the
non-resonance condition~\eqref{eq:non_reso}. Then the function $f$ is a
constant and $c=0$. 
\end{lemma}

\begin{proof}
As $f$ is periodic, we can write it as its Fourier series: denoting $\tau =
(\tau_1,\ldots,\tau_d)$, we have
$$
f(\tau) = \sum_{(k_1,\ldots,k_d) \in \ZZ^d} 
f_k \exp(i k \cdot \tau).
$$
The assumption~\eqref{eq:assump} reads
$$
i \sum_{j=1}^d \sum_{(k_1,\ldots,k_d) \in \ZZ^d} \omega_j 
f_k k_j \exp(i k \cdot \tau) = c.
$$
As $f$ is of class $C^{1+\alpha}$, we have that $| f_k | \leq
C (1+|k|)^{-1-\alpha}$ for some constant $C$, where we recall
$|k| = \sum_{j=1}^d |k_j|.$  Hence, we have that $\dps   
\sum_{k \in \ZZ^d} | f_k k_j | < \infty$, and the above sum is
well-defined.   We deduce that
\begin{equation}
\label{eq:sum_zero}
\forall k \in \ZZ^d, \ \ k \neq 0, \quad
f_k \sum_{j=1}^d \omega_j k_j = 0,
\end{equation}
whereas the identification for $k=0$ leads to $c=0$. 
We infer from~\eqref{eq:sum_zero} and assumption~\eqref{eq:non_reso}
that $f_k = 0$ for any $k \in 
\ZZ^d$, $k \neq 0$. So $f(\tau) = f_0$. This concludes the proof. 
\end{proof}

The equation of order $\eps^0$ reads
\begin{equation}
\label{eq:eps0_pre_multi}
\partial_t S_0 + \sum_{i=1}^d \omega_i \partial_{\tau_i} S_1 
= \frac{P_1^T P_1}{2} + V(q_1 + \partial_{P_1} S_0,0).
\end{equation}
Since $S_0$ does not depend on $\tau$ and $S_1$ is $2 \pi$ periodic with
respect to each $\tau_i$, we can again apply Lemma~\ref{lem:periodic}.
We thus infer from (\ref{eq:eps0_pre_multi}) that
\begin{equation}
\label{eq:S0_pre_multi}
\partial_t S_0 
= \frac{P_1^T P_1}{2} + V(q_1 + \partial_{P_1} S_0,0)
\end{equation}
and $\partial_{\tau_i} S_1 = 0$ for all $1 \leq i \leq d$. 
Equation (\ref{eq:S0_pre_multi}) is supplied with the initial condition
$S_0(t=0,q_1,\alp_2,P_1) = 0$. For each $\alp_2$, we thus recognize
the Hamilton-Jacobi equation for the Hamiltonian function
\begin{equation}
\label{H0_multi}
H_1(q_1,p_1) = \frac{p_1^T p_1}{2} + V(q_1,0).
\end{equation}
So $S_0$ does not depend on $\alp_2$.
In the sequel, we will approximate $S_0(t,q_1,P_1)$ by
\begin{equation}
\label{S0_SE_multi}
S_0^{\rm SE}(t,q_1,P_1) 
= 
S_0(0,q_1,P_1) + t \partial_t S_0(0,q_1,P_1)
=
t \left( \frac{P_1^T P_1}{2} + V(q_1,0) \right),
\end{equation}
which amounts to integrating the Hamiltonian dynamics generated by
(\ref{H0_multi}) with the symplectic Euler algorithm. We have $S_0(t) =
S_0^{\rm SE}(t) + O(t^2)$. 

\medskip

The sequel of the identification is not difficult. The bottom line is
again as follows: since the $\{\tau_i\}$ 
are independent
variables in each $S_k$, we can integrate with respect to each one to
split the equations and close the hierarchy.
Following arguments similar to those presented in~\cite{dcds}, we
find that 
\begin{equation}
\label{eq:S1}
S_1 \equiv 0
\end{equation}
and that $S_2(t) = S_2^{\rm SE}(t) + O(t^2)$, with
\begin{eqnarray}
\nonumber
S_2^{\rm SE}(t,\tau_1,\ldots,\tau_d,q_1,\alp_2,P_1,Y_2)
&=&
\sum_{i=1}^d 
\frac{1}{\omega_i^2} (\nabla_{2,i} V)^T Y_{2,i} \\
\nonumber
&+&
\sum_{i=1}^d \frac{1}{\omega_i^2} 
(\nabla_{2,i} V(q_1+t P_1,0))^T 
\left( (\sin \tau_i) \alp_{2,i} - (\cos \tau_i) Y_{2,i} \right) \\
\nonumber
&-& \frac{t}{2} \sum_{i=1}^d \frac{1}{\omega_i^2} (\nabla_{2,i}
V)^T \nabla_{2,i} V
\\
\label{S2_SE_multi}
&+&
\frac{t}{4} \sum_{i=1}^d \frac{1}{\omega_i^2} 
\left( \alp_{2,i}^T \nabla^2_{2,i} V \alp_{2,i} + Y_{2,i}^T
  \nabla^2_{2,i} V Y_{2,i} \right),
\end{eqnarray}
where the derivatives of $V$ are, unless otherwise mentioned,
evaluated at $(q_1,0)$, and where $\nabla^2_{2,i} V$ is the Hessian
matrix of $V$ with respect to $q_{2,i}$. 

Observe that,
in~\eqref{S2_SE_multi}, there is no term coupling components associated to
different frequencies $\omega_i$ and $\omega_j$, $j \neq i$. This is
reminiscent of the fact that, in the 
ansatz, the fast times $\tau_i$ are independent variables, and that each
$S_k$ is $2 \pi$ periodic with respect to {\em each} $\tau_i$. 

\medskip

Consider now the approximation 
$S_\varepsilon(h) \approx S_\varepsilon^{\rm non-reso}(h)$, with
\begin{equation*}
S_\varepsilon^{\rm non-reso}(h) :=
S_0^{\rm SE}(h) + \eps S_1(h) + \eps^2 S_2^{\rm SE}(h) ,
\end{equation*}
where $S_0^{\rm SE}$, $S_1$ and $S_2^{\rm SE}$ are respectively defined by
(\ref{S0_SE_multi}), (\ref{eq:S1}) and (\ref{S2_SE_multi}).
Using this approximation, we obtain a
symplectic algorithm in the variables $(q_1,x_2,p_1,y_2)$. 
Returning to the original variables $(q_1,q_2,p_1,p_2)$, we obtain the
symplectic Algorithm~\ref{algo:scheme_pre1} outlined below, which we
denote by 
$(Q_1,Q_2,P_1,P_2) = \Psi_h^{\rm non-reso1}(q_1,q_2,p_1,p_2)$.

\bookbox{
\begin{algorithm}[Preconditioned Symplectic Scheme 
$\Psi^{\rm non-reso1}_h(q_1,q_2,p_1,p_2)$]

Set $(q_1,q_2,p_1,p_2) = \left( q_1^n,q_2^n,p_1^n,p_2^n \right)$,
$\tau_i = \omega_i h/\varepsilon$ and 
  perform the following steps:
\begin{enumerate}
\item \text{Change of variables: set} 
$
\quad 
x_2 = q_2, \quad y_2 = p_2, \quad \alp_2 = \Omega x_2/\eps.
$
\item Solve for $(P_1,Y_2)$ in the equations 
$$
\left\{
\begin{array}{rcl}
y_2 
&=& 
\dps
Y_2 
+ \frac{h \varepsilon}{2} \sum_{i=1}^d \frac{1}{\omega_i} \nabla^2_{2,i}
V(q_1,0) \alp_{2,i}  
+
\varepsilon \sum_{i=1}^d \frac{\sin \tau_i}{\omega_i}
\nabla_{2,i} V(q_1+h P_1,0),
\\
p_1 
&=& \dps
P_1 + h \nabla_1 V(q_1,0) 
+ \varepsilon^2 \sum_{i=1}^d 
\frac{1}{\omega_i^2} \nabla_{12i} V(q_1,0) Y_{2,i} 
\\
&+& \dps
\varepsilon^2 \sum_{i=1}^d \frac{1}{\omega_i^2} 
\nabla_{12i} V(q_1+h P_1,0) 
\left( (\sin \tau_i) \alp_{2,i} - (\cos \tau_i) Y_{2,i} \right)  \\
&-& \dps
h \varepsilon^2 \sum_{i=1}^d \frac{1}{\omega_i^2} \nabla_{12i}
V(q_1,0) \nabla_{2,i} V(q_1,0)  
\\
&+& \dps
\frac{h \varepsilon^2}{4} \sum_{i=1}^d \frac{1}{\omega_i^2} 
(\alp_{2,i}^T \nabla_{12i2i} V(q_1,0) \alp_{2,i} 
+ Y_{2,i}^T \nabla_{12i2i} V(q_1,0) Y_{2,i})
.
\end{array}
\right.
$$
\item Set \quad
$\dps
Q_1
=
q_1 
+ h P_1 
+ h \varepsilon^2 \sum_{i=1}^d \frac{1}{\omega_i^2} 
\nabla_{12i} V(q_1+h P_1,0) 
\left( (\sin \tau_i) \alp_{2,i} - (\cos \tau_i) Y_{2,i} \right) .
$
\item Set 
$$
X_2 
=
x_2 + \varepsilon^2 \sum_{i=1}^d \frac{1}{\omega_i^2} \nabla_{2,i} V(q_1,0) 
+ h \varepsilon^2 \sum_{i=1}^d \frac{1}{2 \omega_i^2} 
\nabla^2_{2,i} V(q_1,0) Y_{2,i}
- \varepsilon^2 \sum_{i=1}^d \frac{\cos \tau_i}{\omega_i^2}
\nabla_{2,i} V(q_1+hP_1,0) .
$$
\item Return to the original variables: set $\tau = \Omega h/\eps$ and 
$$
Q_2 = \dps{
(\cos \tau) X_2 +
\varepsilon \Omega^{-1} (\sin \tau) Y_2 },
\quad
P_2 = \dps{
- \frac{\Omega}{\varepsilon} (\sin \tau) X_2 + (\cos \tau) Y_2 }.
$$
\end{enumerate}
\quad Set
$
\left( q_1^{n+1},q_2^{n+1},p_1^{n+1},p_2^{n+1} \right) =
(Q_1,Q_2,P_1,P_2)  .
$
\label{algo:scheme_pre1}
\end{algorithm}
}

Note that, at step 2, we need to solve a system implicit in
$Z = (P_1,Y_2)$, which reads $z = Z + h F(Z) + \eps G(Z)$ with
$z=(p_1,y_2)$. In practice, we use a fixed point method, which converges
in only a few iterations (three iterations in the test-case considered below).

\begin{remark}
High-order derivatives of $V$ appear in Algorithm
\ref{algo:scheme_pre1}. As shown in Section~\ref{sec:finite_diff} above,
they can be replaced by a finite difference approximation in the
generating function.
\end{remark}

Neglecting all terms of order $\eps^3$, the scheme $\Psi^{\rm
non-reso1}_h(q_1,q_2,p_1,p_2)$ is first order in~$h$. A simple, well-known,
manner to get a scheme of higher order is to consider the symmetric form
$$
(Q_1,Q_2,P_1,P_2) = \Psi_h^{\rm non-reso2}(q_1,q_2,p_1,p_2) = 
\left( \Psi_{h/2}^{\rm non-reso1} \right)^* \, \Psi_{h/2}^{\rm non-reso1}
(q_1,q_2,p_1,p_2),
$$
where $\Psi^*$ denotes the adjoint of $\Psi.$
This scheme, denoted Algorithm~\ref{algo:scheme_pre2} in the sequel, is
symplectic, symmetric and, neglecting all terms of order $\eps^3$, of
order 2 in~$h$.

\bookbox{
\begin{algorithm}[Preconditioned Symplectic Scheme 
$\Psi^{\rm non-reso2}_h(q_1,q_2,p_1,p_2)$]

Set $(q_1,q_2,p_1,p_2) = \left( q_1^n,q_2^n,p_1^n,p_2^n \right)$ and 
  perform the following steps:
\begin{enumerate}
\item Set $(\overline{Q}_1, \overline{Q}_2, \overline{P}_1,
  \overline{P}_2) = \Psi_{h/2}^{\rm non-reso1} (q_1,q_2,p_1,p_2)$. 
\item Set $(Q_1,Q_2,P_1,P_2) = \left( \Psi_{h/2}^{\rm non-reso1} \right)^*
(\overline{Q}_1, \overline{Q}_2, \overline{P}_1,\overline{P}_2)$.
\end{enumerate}
\quad Set
$
\left( q_1^{n+1},q_2^{n+1},p_1^{n+1},p_2^{n+1} \right) =
(Q_1,Q_2,P_1,P_2)  .
$
\label{algo:scheme_pre2}
\end{algorithm}
}

\subsection{The non-resonant case: numerical results}
\label{sec:non_reso_num}

We consider a Hamiltonian of the form~\eqref{ham_matrix2}, with $q_1 \in
\RR$ and $q_2 = (q_{2,1},q_{2,2},q_{2,3}) \in \RR^3$, and where
the slow potential energy is
$$
V(q_1,q_2) = 
\left( c + q_{2,1} + q_{2,2} + \gamma q_{2,3} \right)^4
+ 
\frac18 q^2_1 q^2_{2,1} + \frac12 q^2_1
$$
with $c = 1$ and $\gamma = 2.5$. We choose $\Omega = \text{diag} \left(
  1, 1, \sqrt{2} \right)$ as the matrix of fast
frequencies. 

Let $I_j$ denote the energy associated to each fast degree of freedom:
$$
I_j = \frac{p_{2,j}^2}{2} + \frac{\overline{\omega}^2_j \ q_{2,j}^2}{2 \eps^2},
\quad 1 \leq j \leq 3,
$$
with $\overline{\omega}_1 = \overline{\omega}_2 = 1$ and
$\overline{\omega}_3 = \sqrt{2}.$  We first note that
here we use $I_j$ to denote the fast energy, as
opposed to Section~\ref{sec:approach} where it denotes the fast action
(energy divided by frequency).  Of course, for the constant frequency
case here, these two quantities differ only by a multiplicative constant.
In addition, we use here a different convention than that of
Section~\ref{sec:non_reso} on how the eigenvalues of $\Omega$ are
numbered.  It is well-known (see~\cite[Sec. XIII.9]{HLW}) that the
quantities
\begin{equation}
\label{eq:I_multi}
I = \sum_{j=1}^3 I_j \quad \text{and} \quad I_3
\end{equation}
are adiabatic invariants of the dynamics. In contrast, $I_1$ and $I_2$,
associated to the same frequency, are {\em not} adiabatic invariants, 
although their sum, $I_1+I_2 = I - I_3$, is. 

\medskip

\begin{figure}[htbp]
\centerline{
\input{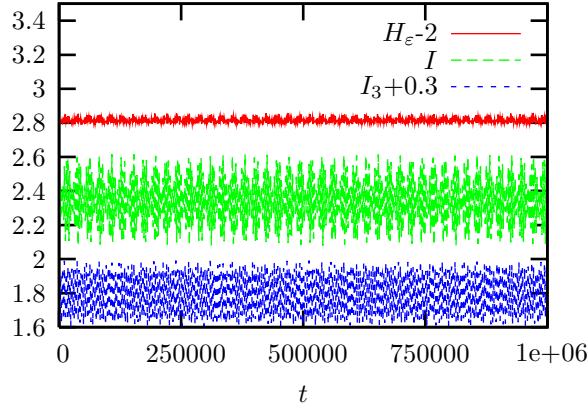}
}
\caption{\small
Energy and adiabatic invariants $I$ and $I_3$ (for convenience, we plot
$H_\eps-2$, $I$ and $I_3+0.3$) along the trajectory computed with
Algorithm~\ref{algo:scheme_pre2} ($\eps = 1/70$ and $h=10\eps$).}
\label{fig:drift_non_reso}
\end{figure}

We first choose $\eps=1/70$ and $h=10 \eps$, and we monitor the evolution
of the energy and adiabatic invariants up to time $T=10^6$ on the
numerical trajectory computed with Algorithm~\ref{algo:scheme_pre2}.
Results are shown in Fig.~\ref{fig:drift_non_reso}. 
We observe no drift.

\begin{figure}[htbp]
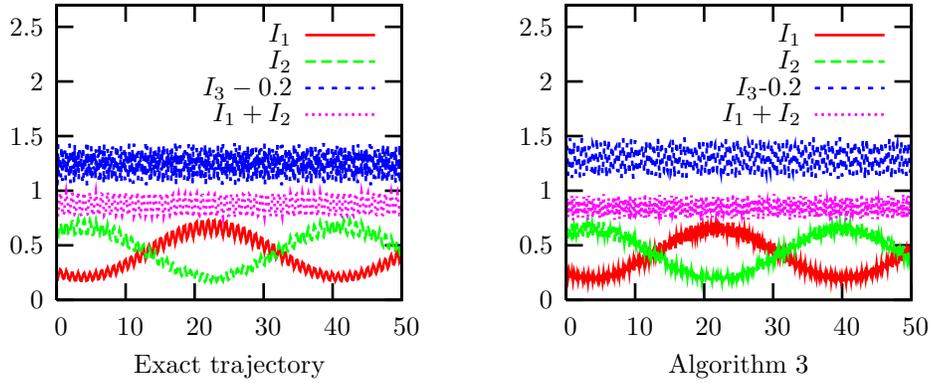

\centerline{
\begin{tabular}{cc}
\input{figs/echange_exact_non_reso.tex} & 
\input{figs/echange_mts_non_reso.tex} 
\end{tabular}
}
\caption{\small
Preservation of the adiabatic invariants $I_1+I_2$ and $I_3$, and exchange
between $I_1$ and $I_2$, for $\eps = 1/70$ (Algorithm
\ref{algo:scheme_pre2} has been used with $h=10 \eps$).}  
\label{fig:echange_non_reso}
\end{figure}

For the same parameters, we show in Fig.~\ref{fig:echange_non_reso} the
evolution of $I_j$ over the time window $[0,50]$. As expected, $I_3$ is
preserved, as well as $I_1+I_2$. We observe that Algorithm
\ref{algo:scheme_pre2} correctly reproduces the exchange between $I_1$
and $I_2$. 

\medskip

\begin{figure}[htbp]
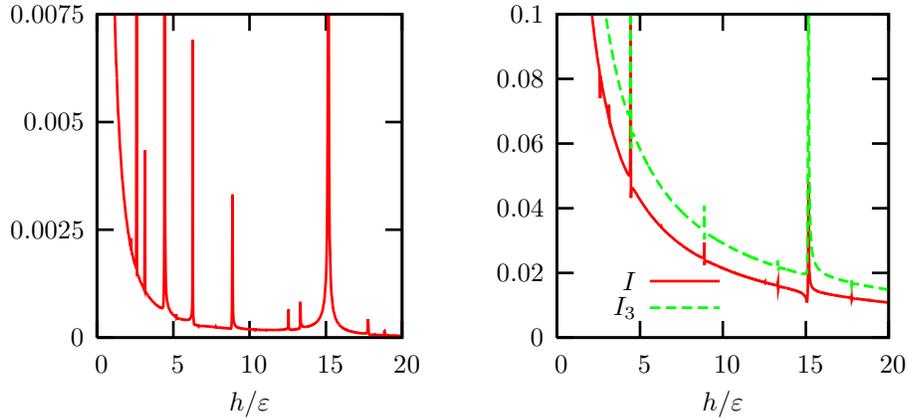

\centerline{
\input{figs/comparison_h0.02_non_reso_energy.tex}
\input{figs/comparison_h0.02_non_reso_inv_adiab.tex} 
}
\caption{\small 
Maximum variations (\ref{eq:estim}) of the energy (left) and of the
adiabatic invariants $I$ and $I_3$ (right) on the time interval
$[0,10^4]$, for several $\eps$ ($h = 0.02$), for Algorithm
\ref{algo:scheme_pre2}. 
\label{fig:h_fixed_eps_vary_non_reso}}
\end{figure}

We now focus on the robustness of Algorithm~\ref{algo:scheme_pre2} as
$\eps$ decreases. We
set the time step to $h=0.02$ and consider the variations of the energy
and of the adiabatic invariants~\eqref{eq:I_multi},
\begin{equation}
\label{eq:estim}
\max_{t \in [0,10^4]} \frac{\left| H(t) - H(0)
  \right|}{H(0)},
\quad
\max_{t \in [0,10^4]} \frac{\left| I(t) - I(0)
  \right|}{I(0)},
\quad
\max_{t \in [0,10^4]} \frac{\left| I_3(t) - I_3(0)
  \right|}{I_3(0)},
\end{equation}
over the time interval $t \in [0,10^4]$, for
stiffness $\eps$ varying between $10^{-3}$ to 1. Results are shown on
Fig.~\ref{fig:h_fixed_eps_vary_non_reso}. Only a few resonances can be
seen, and they are extremely peaked. In addition, the algorithm
is very stable as $\eps$ decreases to 0.

\subsection{The resonant case: derivation of the scheme}
\label{sec:reso}

In this section, we again consider the Hamiltonian~\eqref{ham_matrix2}, with a
diagonal matrix $\Omega \in \RR^{f \times f}$. In contrast to
Section~\ref{sec:non_reso}, we assume here that some entries of the diagonal
of $\Omega$ are resonant. To simplify the notation and
the analysis, we assume that $\Omega$ is of the
form~\eqref{eq:def_Omega_diago} with just $d=2$, namely 
$$
\Omega = \left(
\begin{array}{cc}
\omega_a \Id_{m_a} & 0 
\\
0 & \omega_b \Id_{m_b}
\end{array}
\right),
$$
with $m_a + m_b = f$, $0 < \omega_a < \omega_b$, and 
\begin{equation}
\label{eq:depend}
\omega_b = \frac{\beta}{\alpha} \omega_a, \ \ \text{with} \
\alpha \in \NN^\star, \ \beta \in \NN^\star, 
\ \alpha \neq \beta, \ \text{and $\alpha$ and $\beta$ coprime}.
\end{equation}
Likewise, we write $q_2=(q_{2,a},q_{2,b})$ with $q_{2,a} \in \RR^{m_a}$ and
$q_{2,b} \in \RR^{m_b}$, and $p_2=(p_{2,a},p_{2,b})$.

\medskip

As in Section~\ref{sec:non_reso}, we first consider the time-dependent change
of variables~\eqref{eq:cov_multi}, which we write
$(x_2,y_2) = \chi(t,q_2,p_2) \in \RR^{2f}$. 
The dynamics on $(q_1,x_2,p_1,y_2)$ is a Hamiltonian
dynamics with the time-dependent Hamiltonian
$$
H_\varepsilon^{\rm reso}(t,q_1,x_2,p_1,y_2) = 
\frac{p_1^T p_1}{2} + W_\eps^{\rm reso} \left(\frac{\omega_a
  t}{\alpha \eps},q_1,x_2,y_2 \right), 
$$
where 
$$
W_\eps^{\rm reso}(\tau,q_1,x_2,y_2)
= V \left[q_1,
(\cos \alpha \tau) x_{2,a} 
+ \frac{\varepsilon}{\omega_a} (\sin \alpha \tau) y_{2,a},
(\cos \beta \tau) x_{2,b} 
+ \frac{\varepsilon}{\omega_b} (\sin \beta \tau) y_{2,b}
\right].
$$
We let $\overline{S}_\eps(t,q_1,x_2,P_1,Y_2)$ denote the solution to the
Hamilton-Jacobi equation associated with $H_\varepsilon^{\rm reso}$ and
perform the change of variables and of unknown function
$$
\alp_2 = \frac{\Omega}{\eps} x_2 \quad \text{and} \quad 
S_\eps(t,q_1,\alp_2,P_1,Y_2) = 
\overline{S}_\eps\left(t,q_1,\eps \Omega^{-1} \alp_2,P_1,Y_2\right).
$$
The function $S_\eps$ satisfies
\begin{equation}
\label{hj3_reso}
\left\{
\begin{array}{l}
\dps
\partial_t S_\eps = \frac{P_1^T P_1}{2} +
W_\eps^{\rm reso} \left(\frac{\omega_a t}{\alpha \eps},
q_1 + \partial_{P_1} S_\eps,
\eps \Omega^{-1} \alp_2 + \partial_{Y_2} S_\eps, Y_2 \right),
\\ \noalign{\vskip 4pt}
S_\eps(0,q_1,\alp_2,P_1,Y_2) = 0.
\end{array}
\right.
\end{equation}
We make the ansatz
\begin{eqnarray}
\label{ansatz2_reso}
S_\eps(t,q_1,\alp_2,P_1,Y_2) &=& 
S_0 \left( t,\tau,q_1,\alp_2,P_1,Y_2 \right) 
+
\eps S_1 \left( t,\tau,q_1,\alp_2,P_1,Y_2 \right) 
\\
\nonumber
&+& \text{higher order terms in $\eps^k$, $k \geq 2$},
\end{eqnarray}
where the fast time $\tau$ is defined by
\begin{equation}
\label{eq:tau_reso}
\tau = \frac{t \omega_a}{\alpha \eps} = \frac{t \omega_b}{\beta \eps},
\end{equation}
and where the functions $(S_k)_{k \geq 0}$ are assumed to be $2 \pi$ 
periodic with respect to $\tau$. 

\begin{remark}
In the case where $V$ does not depend on $q_2$, the solution
to~\eqref{hj3_reso} can be analytically identified and is indeed of
the form~\eqref{ansatz2_reso}. 
\end{remark}

We now insert (\ref{ansatz2_reso}) in (\ref{hj3_reso}), identify the first
variable of $W_\eps^{\rm reso}$ with the fast time $\tau$, and expand in
powers of $\eps$. 
As in Section~\ref{sec:non_reso}, $S_0$ is independent from $Y_2$ and
$\tau$. The equation of order $\eps^0$ reads
\begin{equation}
\label{eq:eps0_pre_reso}
\partial_t S_0 + \frac{\omega_a}{\alpha} \partial_\tau S_1 
= \frac{P_1^T P_1}{2} + V(q_1 + \partial_{P_1} S_0,0).
\end{equation}
Since $S_0$ does not depend on $\tau$ and $S_1$ is $2 \pi$ periodic with
respect to $\tau$, we infer from (\ref{eq:eps0_pre_reso}) that
\begin{equation}
\label{eq:S0_pre_reso}
\partial_t S_0 
= \frac{P_1^T P_1}{2} + V(q_1 + \partial_{P_1} S_0,0)
\end{equation}
and $\partial_\tau S_1 = 0$.
Equation (\ref{eq:S0_pre_reso}) is supplied with the initial condition
$S_0(t=0,q_1,\alp_2,P_1) = 0$. We again recognize the Hamilton-Jacobi
equation for the Hamiltonian function~\eqref{H0_multi}. In the sequel,
$S_0$ is thus again approximated by~\eqref{S0_SE_multi}. 

\medskip

Following the same arguments as in~\cite{dcds}, we next proceed with the
sequel of the identification, using the fact that, since $\alpha \neq
\beta$ in~\eqref{eq:depend},
$$
\int_0^{2 \pi} \cos \alpha \tau \ \cos \beta \tau \ d\tau= 0.
$$
As a consequence, at least at the orders of the expansion that we
consider, no coupling appears between the terms associated to the
frequency $\omega_a$ and those associated with the frequency
$\omega_b$. We find that
\begin{equation}
\label{eq:S1_reso}
S_1 \equiv 0
\end{equation}
and that $S_2(t) = S_2^{\rm SE}(t) + O(t^2)$, with
\begin{eqnarray}
\nonumber
S_2^{\rm SE}(t,\tau,q_1,\alp_2,P_1,Y_2) 
&=&
\frac{1}{\omega_a^2} (\nabla_{2,a} V)^T Y_{2,a}
+ \frac{1}{\omega_b^2} (\nabla_{2,b} V)^T Y_{2,b}
\\
\nonumber
&+&
\frac{1}{\omega_a^2} 
(\nabla_{2,a} V(q_1+t P_1,0))^T 
\left( (\sin \alpha \tau) \alp_{2,a} - (\cos \alpha \tau) Y_{2,a}
\right)
\\
\label{S2_SE_reso}
&+&
\frac{1}{\omega_b^2} 
(\nabla_{2,b} V(q_1+t P_1,0))^T 
\left( (\sin \beta \tau) \alp_{2,b} - (\cos \beta \tau) Y_{2,b} \right) 
\\ 
\nonumber
&-&
\frac{t}{2} \left( \frac{1}{\omega_a^2} (\nabla_{2,a}
V)^T \nabla_{2,a} V + \frac{1}{\omega_b^2} (\nabla_{2,b}
V)^T \nabla_{2,b} V \right)
\\
\nonumber
&+&
\frac{t}{4 \omega_a^2} 
\left( \alp_{2,a}^T \nabla^2_{2,a} V \alp_{2,a} + Y_{2,a}^T
  \nabla^2_{2,a} V Y_{2,a} \right)
\\
\nonumber
&+&
\frac{t}{4\omega_b^2} 
\left( \alp_{2,b}^T \nabla^2_{2,b} V \alp_{2,b} + Y_{2,b}^T
  \nabla^2_{2,b} V Y_{2,b} \right),
\end{eqnarray}
where the derivatives of $V$ are evaluated at $(q_1,0)$ unless
otherwise mentioned, and where $\nabla^2_{2,a} V$ is the Hessian
matrix of $V$ with respect to $q_{2,a}$. 

\medskip

Observe that,
in~\eqref{S2_SE_reso}, there is no term coupling components associated to
different frequencies. This is reminiscent of the fact that, at the
$O(\eps^2)$ term that we consider here, coupling terms can only come
from terms containing products of the form $\cos \alpha \tau \ \cos
\beta \tau$, the average of which vanishes. In contrast, to identify $S_3$, we need to handle 
terms of the form $(\cos \alpha \tau)^s \ (\cos \beta \tau)^r$ (with $r+s
= 3$), whose average may not vanish (e.g. in the case $\alpha = 1$,
$\beta = 2$, $s=2$ and $r=1$). 

\bigskip

We next observe that the generating functions $S_0^{\rm SE}$, $S_1$ and
$S_2^{\rm SE}$ that we have identified in this resonant case, defined by  
(\ref{S0_SE_multi}), (\ref{eq:S1_reso}) and (\ref{S2_SE_reso}), where the
fast time $\tau$ is given by~\eqref{eq:tau_reso}, are equal to the generating
functions $S_0^{\rm SE}$, $S_1$ and $S_2^{\rm SE}$ identified in the
non-resonant case, see~(\ref{S0_SE_multi}), (\ref{eq:S1}) and
(\ref{S2_SE_multi}), where the fast times $\tau_i$ are given
by~\eqref{eq:tau_multi}. As a consequence, even though the frequencies
$\omega_a$ and $\omega_b$ 
are here resonant, we again obtain the algorithms~\ref{algo:scheme_pre1}
and~\ref{algo:scheme_pre2} proposed above.

\subsection{The resonant case: numerical results}
\label{sec:multi_gene_num}

In Sections~\ref{sec:non_reso} and~\ref{sec:reso}, we have derived
algorithms to integrate the dynamics in the non-resonant case and in the
resonant case, respectively. As explained at the end of the previous
section, when expanding up to $O(\eps^2)$, these two cases turn
out to yield the same algorithm. A natural intuition, which is
confirmed by a careful identification similar to the ones performed
above, is that this algorithm is also appropriate in the case when
the matrix of fast frequencies is of the form~\eqref{eq:def_Omega_diago},
where the $\omega_i$ again satisfy $0 < \omega_1 < \omega_2 < \ldots <
\omega_d$, but where, in contrast to the situation considered in
Sections~\ref{sec:non_reso} and~\ref{sec:reso}, they are now
possibly resonant, and $d$ is arbitrary.

Consider a system of the form~\eqref{ham_matrix2}, with $q_1 \in
\RR$ and $q_2 = (q_{2,1},q_{2,2},q_{2,3},q_{2,4}) \in \RR^4$, and where
the slow potential energy is
$$
V(q_1,q_2) = 
\left( c + q_{2,1} + q_{2,2} + q_{2,3} + \gamma q_{2,4} \right)^4
+ 
\frac18 q^2_1 q^2_{2,1} + \frac12 q^2_1
$$
with $c = 1$ and $\gamma = 2.5$. We choose the matrix of fast
frequencies to be $\Omega = \text{diag} \left( 1, 1, \sqrt{2}, 2
\right)$.  A very similar test-case 
has already been studied in~\cite{rennes}
and~\cite[Sec. XIII.9.1]{HLW}, where the chosen value of $c=1/20$
causes the exchange between $I_1$ and $I_2$ to occur on a slower time scale. 
 It obviously enters the framework of the
current section, but not the framework of Section~\ref{sec:non_reso}, as
the first and the last frequencies are resonant. In this
case, the adiabatic invariants are
\begin{equation}
\label{eq:I_reso}
I = \sum_{j=1}^4 I_j \quad \text{and} \quad I_3,
\end{equation}
where $I_j$ is the energy associated to the $j$th fast degree of freedom:
$$
I_j = \frac{p_{2,j}^2}{2} + \frac{\overline{\omega}^2_j \ q_{2,j}^2}{2 \eps^2},
\quad 1 \leq j \leq 4,
$$
with $\overline{\omega}_1 = \overline{\omega}_2 = 1$,
$\overline{\omega}_3 = \sqrt{2}$ and $\overline{\omega}_4 = 2$. 

\medskip

We first choose $\eps=1/70$ and $h=10 \eps$, and we monitor the evolution
of the energy and adiabatic invariants up to time $T=10^6$ on the
numerical trajectory computed with Algorithm~\ref{algo:scheme_pre2}.
Results are very similar to those shown on
Fig.~\ref{fig:drift_non_reso}: we do not observe any drift.

\medskip

\begin{figure}[htbp]
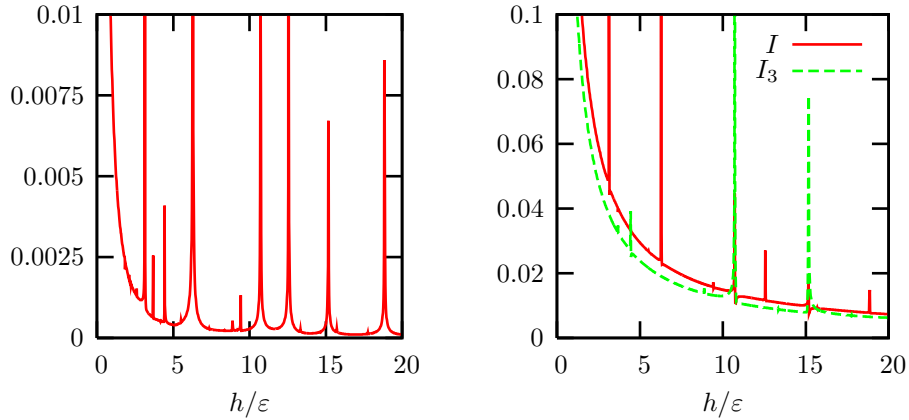

\centerline{
\input{figs/comparison_h0.02_multi_energy.tex}
\input{figs/comparison_h0.02_multi_inv_adiab.tex} 
}
\caption{\small 
Maximum variations (\ref{eq:estim}) of the energy (left) and of the 
adiabatic invariants $I$ and $I_3$ (right) on the
time interval $[0,10^4]$, for several $\eps$ ($h = 0.02$), for Algorithm
\ref{algo:scheme_pre2}. There is an increase in resonances compared
to Fig.~\ref{fig:h_fixed_eps_vary_non_reso}, though as before they are very
tightly peaked.} 
\label{fig:h_fixed_eps_vary_multi}
\end{figure}

We now focus on the robustness of Algorithm~\ref{algo:scheme_pre2}, as
$\eps$ decreases. We
set the time step to $h=0.02$, and consider the
variations~\eqref{eq:estim} of the energy and of the adiabatic
invariants~\eqref{eq:I_reso},  
over the time interval $t \in [0,10^4]$, for
stiffness $\eps$ varying between $10^{-3}$ to 1. Results are shown on
Fig.~\ref{fig:h_fixed_eps_vary_multi}. When $\eps$ decreases to 0, the
algorithm performs better and better, except for a few peaked
resonances.

\medskip

We conclude this section by discussing the exchange of fast energies.
In the case at hand, it is well-known (see~\cite[Sec. XIII.9]{HLW})
that some exchange occurs between $I_1$, $I_2$ (which are both
associated to the same frequency $\overline{\omega}_1$) and $I_4$
(which is associated to a frequency $\overline{\omega}_4$ resonant
with $\overline{\omega}_1$).  The exchange between $I_1$ and $I_2$ is
$O(1)$ over timescales of $O(\eps^{-1})$, which is stronger than the
exchange between $I_4$ and $I_1 + I_2$, which is only $O(\eps)$ over
timescales of $O(\eps^{-1})$.  For some choices of parameters, there
is an $O(1)$ exchange between $I_4$ and $I_1 + I_2$ that occurs on the
long time scale $O(\eps^{-2}).$  It turns out that, on the numerical
trajectory computed using Algorithm~\ref{algo:scheme_pre2}, the
exchange between $I_1 + I_2$ and $I_4$ on long time scales is not
reproduced, and $I_4$ is almost preserved.  This is due to the fact
that (i) resonant frequencies are handled by the algorithm we derived
as if they were non-resonant (see Section~\ref{sec:reso}) and (ii) in
the non-resonant case, no exchange occurs between the fast energies.
If we were to expand the generating function to third order in $\eps,$
certain exchange terms among the fast terms would be included in the
resulting algorithm, which could potentially capture the exchange
between $I_1 + I_2$ and $I_4.$

\begin{figure}[htbp]
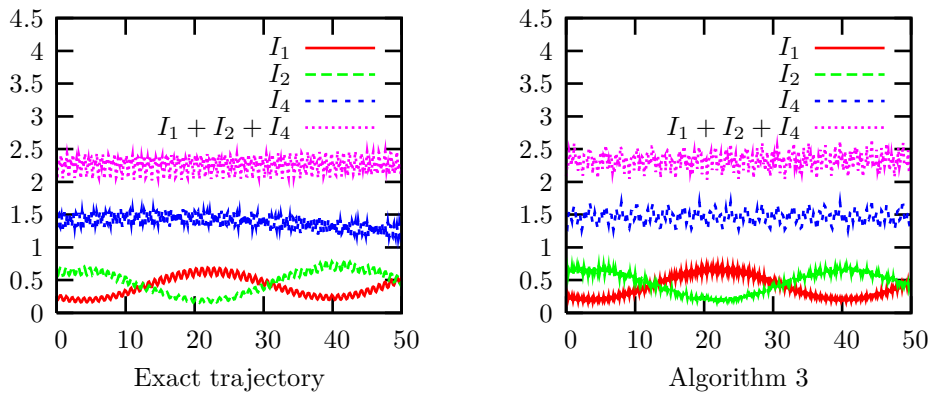

\centerline{
\begin{tabular}{cc}
\input{figs/echange_exact_reso.tex} & 
\input{figs/echange_mts_reso.tex} 
\end{tabular}
}
\caption{\small
The exchange between $I_1$ and $I_2$ on the time scale $O(\eps^{-1})$ is
accurately reproduced by Algorithm~\ref{algo:scheme_pre2}, as well as
the preservation of the adiabatic invariant $I_1+I_2+I_4$.  The
$O(\eps)$ exchange between $I_4$ and $I_1 + I_2$ is not captured. (We
take $\eps=1/70$, and we use $h=10 \eps$ for
Algorithm~\ref{algo:scheme_pre2}).}  
\label{fig:echange_multi}
\end{figure}

Yet, the fact that we chose not to include these terms did not
destroy the preservation of energy and
adiabatic invariants, which are accurately recovered, as pointed out
above.  In addition, the exchange between $I_1$ and $I_2$, that occurs
on the time scale $O(\eps^{-1})$, is accurately reproduced by
Algorithm~\ref{algo:scheme_pre2}, as shown on
Fig.~\ref{fig:echange_multi}.  

\section{The extensible pendulum}
\label{sec:pendulum}

In the previous sections, we have considered Hamiltonians of the
form~\eqref{ham_general}, where the leading order behavior of the fast
variables is that of a harmonic oscillator, with either a slowly varying
scalar frequency (in Section~\ref{sec:approach}), or a constant matrix of
fast frequencies (in Section~\ref{sec:multi}).  On the other hand, 
many examples of interest are not harmonic, and this technique cannot 
handle completely general fast forces. 
An intermediate class of system is the one discussed
in~\cite[Sec. XIV.3]{HLW}, where the Hamiltonian reads
$$
H(q,p) = \frac12 p^T p + V_{\rm slow}(q) + \frac{1}{\eps^2} U_{\rm fast}(q)
$$
with fast potential $U_{\rm fast}$ that has a
change of variables 
$y = (y_1,y_2)= \chi(q)$ so that 
$$
U_{\rm fast}(q) = \frac12 y_2^T \Omega(y_1)^2 y_2,
$$ 
where $\Omega(y_1)$ is a $f \times f$ matrix. So, up to a change of
variables, the fast potential 
energy is of the form considered in~\eqref{ham_general}. However, the
map $(q,p) \mapsto (y,p)$ is in general not
symplectic. Hence we cannot work with the variables $(y,p)$, as the
dynamics in these variables is not Hamiltonian. We thus need
to also change variables for the momenta, so that the map
$(q,p) \mapsto (y,v)$, where $v$ are the new momenta, is
symplectic, and the dynamics in $(y,v)$ is Hamiltonian. A standard
consequence of this change of variables for the momenta is that the
kinetic energy turns out to depend on the positions $y$. The Hamiltonian
reads
\begin{equation}
\label{eq:complique}
H(y,v) = T(y,v) + V_{\rm slow}(\chi^{-1}(y)) + 
\frac{1}{2\eps^2} y_2^T \Omega(y_1)^2 y_2.
\end{equation}
Although the
fast potential is harmonic, the fast Hamiltonian is not necessarily
close to that of a harmonic oscillator.

\smallskip

In this section, we consider a particular case
(see~\eqref{eq:ham_pendule} below) that cannot be written in
the form~\eqref{ham_general}. After a global change of variables, in
positions and momenta, we transform the Hamiltonian to the
form~\eqref{eq:complique} (see~\eqref{eq:ham_pendule_int} below), where
the kinetic energy turns out to still 
be a quadratic form of the momenta: $\dps T(y,v) = \frac12 v^T M^{-1}(y)
v$. The resulting mass matrix $M(y)$ depends on the position. 
We will show that our strategy can still handle this case and yields an
efficient algorithm.

\subsection{Transforming to internal coordinates}
\begin{figure}[htbp]
\psfrag{a}{$a$}
\psfrag{r}{$\overline{r}$}
\psfrag{qx}{$(q_x,q_y)$}
\centerline{
\includegraphics[width=2cm]{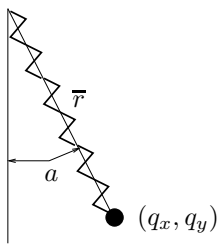}
}
\caption{Extensible pendulum test-case: a particle at position
  $(q_x,q_y)$ is attached to the origin using an extensible spring. We
  denote by $\overline{r}$ the spring length and by $a$ the angle between the
  pendulum and the vertical.} 
\label{fig:pendule_image}
\end{figure}

We consider an extensible pendulum in two
dimensions (see 
Fig.~\ref{fig:pendule_image}). The potential energy of the system is the
sum of two terms, the spring energy and a term $W(a)$ that
depends on the angle $a$. Denote by $q_x$
and $q_y$ the Euclidean coordinates of the particle and by $p_x$ and $p_y$ the
corresponding momenta. The Hamiltonian of the system reads
\begin{equation}
\label{eq:ham_pendule}
H_{\rm cartesian}(q_x,q_y,p_x,p_y) = \frac12 (p_x^2 + p_y^2) + 
\frac{1}{2 \eps^2} (\overline{r} - r_0)^2 + W(a),
\end{equation}
where $\overline{r} = \sqrt{q_x^2 + q_y^2}$ is the length of the spring,
$r_0$ is the equilibrium length of the spring, $\eps$ is a small
parameter, and $W$ is a $2\pi$-periodic non-negative function.

For illustration, this example can be considered as an oversimplified
model of a molecular system, in which the spring-like potential 
$(2 \eps^2)^{-1} (\overline{r} - r_0)^2$ models each bond length between
neighboring atoms, and with some potential associated with each
bond angle. 

The fast oscillating term $(2 \eps^2)^{-1} (\overline{r} - r_0)^2$
of the Hamiltonian~\eqref{eq:ham_pendule} is not harmonic with respect
to the {\em cartesian} degrees of freedom $(q_x,q_y)$, and thus the 
strategy described in Section~\ref{sec:approach} cannot be directly  
applied to~\eqref{eq:ham_pendule}. However, as is often the case in
molecular 
dynamics models, the fast oscillating term of~\eqref{eq:ham_pendule} is
harmonic in the {\em internal} degrees of freedom, the bond
length $\overline{r}$ and the bond angle $a$. It is thus natural to introduce
polar coordinates $(a,\overline{r})$, with 
$$
q_x = \overline{r} \cos a \quad \text{and} \quad q_y = \overline{r} \sin a,
$$
so that the fast oscillating term of the potential energy is {\em
  harmonic} in these coordinates. 
We now introduce momenta $(p_a,p_r)$ associated to
$(a,\overline{r})$, such that the map $(q_x,q_y,p_x,p_y) \mapsto
(a,\overline{r},p_a,p_r)$ is symplectic. Following~\cite[Example~VI.5.2]{HLW}, this leads to choosing
$$
p_r = p_x \cos a + p_y \sin a 
\quad \text{and} \quad
p_a = - p_x \overline{r} \sin a + p_y \overline{r} \cos a.
$$
In addition, rather than working with
$\overline{r}$, we work in the sequel with
$$
r = \overline{r} - r_0,
$$
so that the fast oscillating harmonic term attains its minimum at $r=0$,
as in~\eqref{ham_general}. Without loss of generality, we furthermore
assume that $r_0 = 1$. The Hamiltonian then reads
\begin{equation}
\label{eq:ham_pendule_int}
H(a,r,p_a,p_r)
=
\frac{p_r^2}{2} + \frac{p_a^2}{2 (1 + r)^2}
+ 
\frac{1}{2 \eps^2} r^2 + W(a).
\end{equation}
Although the mass matrix depends on the fast position $r$, we will show
that our general strategy, originally developed for~\eqref{ham_general},
can still be used. 

\begin{remark}
As above, we choose initial conditions depending on $\eps$ such that the
energy is bounded. As $W$ is bounded from below, this implies that, at any
time $t$, $r(t) = O(\eps)$. As a consequence, one could
approximate~\eqref{eq:ham_pendule_int} by
$$
H_{\rm decoup}(a,r,p_a,p_r)
=
\frac{p_r^2}{2} + \frac{p_a^2}{2}
+ 
\frac{1}{2 \eps^2} r^2 + W(a).
$$
In this case, $(r,p_r)$ and $(a,p_a)$ are 
decoupled. In the following, we do not make this approximation, and we
work with~\eqref{eq:ham_pendule_int}. 
\end{remark}

\subsection{Derivation of the symplectic scheme}

Starting from~\eqref{eq:ham_pendule_int}, we now follow our usual
strategy: we first precondition the fast variables, and next consider
the Hamilton-Jacobi form of the equations. Observe that, in the case at hand, the fast frequency is
constant. We thus follow the approach described in~\cite{dcds} and
consider the time-dependent change of variable $(r,p_r) \mapsto
(\overline{b},p_b)$ defined by
\begin{equation}
\label{eq:cov}
r = \overline{b} \cos \frac{t}{\eps} + \eps p_b \sin \frac{t}{\eps}
\quad \text{and} \quad
p_r = - \frac{\overline{b}}{\eps} \sin \frac{t}{\eps} + 
p_b \cos \frac{t}{\eps}.
\end{equation}
In the variables $(a,\overline{b},p_a,p_b)$, the dynamics is again a
Hamiltonian dynamics for the Hamiltonian $\dps
H_{\rm pend}\left(\frac{t}{\eps},a,b,p_a,p_b\right)$, with
$$
H_{\rm pend}(\tau,a,\overline{b},p_a,p_b) = 
\frac{p_a^2}{2 (1 + \overline{b} \cos \tau + \eps p_b \sin \tau)^2} + W(a).
$$

Let $\overline{S}_\eps(t,a,\overline{b},P_a,P_b)$ solve the
Hamilton-Jacobi equation associated to $H_{\rm pend}$. 
As $r$ and $\overline{b}$ are of order $\eps$, we make the change of
variables and of unknown function:
$$
b = \frac{\overline{b}}{\eps} \quad \text{and} \quad 
S_\eps(t,a,b,P_a,P_b) = 
\overline{S}_\eps\left(t,a,\eps b,P_a,P_b\right),
$$
so that $S_\eps$ satisfies
\begin{equation}
\label{hj3}
\partial_t S_\eps = H_{\rm pend}\left(\frac{t}{\eps},a + \partial_{P_a} S_\eps,
\eps b + \partial_{P_b} S_\eps,P_a,P_b\right),
\quad 
S_\eps(0,a,b,P_a,P_b) = 0.
\end{equation}
We make the ansatz
\begin{eqnarray}
\label{ansatz2}
S_\eps(t,a,b,P_a,P_b) &=& 
S_0 \left(t,\tau,a,b,P_a,P_b \right) +
\eps S_1 \left(t,\tau,a,b,P_a,P_b \right) 
\\
\nonumber
&+& \text{higher order terms in $\eps^k$, $k \geq 2$},
\end{eqnarray}
where the fast time $\tau$ is defined by
$$
\tau  = \frac{t}{\eps},
$$
and where the functions $(S_k)_{k \geq 0}$ are assumed to be $2 \pi$ 
periodic in $\tau$. We now insert (\ref{ansatz2}) in (\ref{hj3}), identify the first
variable of $H_{\rm pend}$ with the fast time $\tau$, and expand in powers of
$\eps$. 

\medskip

The identification follows the same lines as in the previous
sections. We obtain that $S_0$ is independent of $b$, $P_b$ and $\tau$,
and satisfies $S_0(t) = S_0^{\rm SE}(t) + O(t^2)$ with 
\begin{equation}
\label{S0_SE}
S_0^{\rm SE}(t,a,P_a) 
=
t \left( \frac{P_a^2}{2} + W(a) \right),
\end{equation}
while 
\begin{equation}
\label{eq:S1_pend}
S_1 \equiv 0.
\end{equation}
Using the method of characteristics, we obtain that 
$S_2(t) = S_2^{\rm SE}(t) + O(t^2)$, with
\begin{eqnarray}
\label{S2_SE}
S_2^{\rm SE}(t,\tau,a,b,P_a,P_b) &=&
P_a^2 \left[ P_b \cos \tau - b \sin \tau \right]
- P_b (P_a + t W'(a))^2 
\\
\nonumber
&+& t \left[
\frac34 (b^2 + P_b^2) (P_a + t W'(a))^2 - \frac12 (P_a + t W'(a))^4
\right].
\end{eqnarray}
Collecting~\eqref{S0_SE}, \eqref{eq:S1_pend} and~\eqref{S2_SE}, we obtain
the following approximation of the generating function:
$$
S_\eps(h,a,b,P_a,P_b) \approx
\widetilde{S_\eps}(h,a,b,P_a,P_b)
=
S_0^{\rm SE}(h,a,P_a) + \eps^2 S_2^{\rm SE}(h,\frac{h}{\eps},a,b,P_a,P_b).
$$
Returning to the variables $(a,r,p_a,p_r)$ of the  
Hamiltonian in~\eqref{eq:ham_pendule_int}, we obtain a symplectic scheme,
called Algorithm~\ref{algo:pendulum} in the sequel, which we denote by
$(a^{n+1},r^{n+1},p_a^{n+1},p_r^{n+1}) = 
\Psi_h^{\rm pendulum}(a^{n},r^{n},p_a^{n},p_r^{n})$. 

\bookbox{
\begin{algorithm}[Symplectic scheme $\Psi^{\rm pendulum}_h(a,r,p_a,p_r)$]

Set $(a,r,p_a,p_r) = \left( a^n,r^n,p_a^n,p_r^n \right)$,
$\tau = h/\varepsilon$ and perform the following steps:
\begin{enumerate}
\item \text{Change of variables: set}  
$ 
\quad b = r/\eps, \quad p_b = p_r.
$
\item Solve for $(P_a,P_b)$ in the equations 
$$
\left\{
\begin{array}{rcl}
p_b 
&=& 
\dps{
P_b 
- \varepsilon P_a^2 \sin \tau 
+ \frac32 h \varepsilon b \, (P_a + h W'(a))^2,
}
\\ \noalign{\vskip 4pt}
p_a 
&=& 
\dps{
P_a + h W'(a) 
- 2 h \eps^2 P_b \, (P_a + h W'(a)) W''(a) }
\\ \noalign{\vskip 4pt}
&+& \dps{
h^2 \eps^2 \left[ \frac32 (b^2+P_b^2) (P_a + h W'(a)) - 2
  (P_a + h W'(a))^3 \right] W''(a) .
}
\end{array}
\right.
$$
\item Set
\begin{eqnarray*}
A &=& a + h P_a + 2 \eps^2 P_a (P_b \cos \tau - b \sin \tau) - 
2 \eps^2 P_b (P_a + h W'(a))
\\
&+& 
h \varepsilon^2 \left[ \frac32 (b^2+P_b^2) (P_a + h W'(a)) - 
2 (P_a + h W'(a))^3 
\right].
\end{eqnarray*}
\item Set 
$\dps 
B =
b + \varepsilon P_a^2 \cos \tau - \eps (P_a + h W'(a))^2 
+ \frac32 h \varepsilon P_b (P_a + h W'(a))^2 .
$
\item Return to the original variables: \ \
$
R = \eps B \cos \tau +
\varepsilon P_b \sin \tau ,
\quad
P_r =
- B \sin \tau  + P_b \cos \tau .
$
\end{enumerate}
\quad Set
$
\left( a^{n+1},r^{n+1},p_a^{n+1},p_r^{n+1} \right) =
(A,R,P_a,P_r)  .
$
\label{algo:pendulum}
\end{algorithm}
}
To solve the implicit equations for $(P_a,P_b)$, we first observe that,
once $P_a$ is known, $P_b$ can be explicitly determined. We thus first
recast the problem as a scalar implicit equation on $P_a$ and solve it using a
fixed point algorithm before determining $P_b$. 

Neglecting all terms of order $\eps^3$, the scheme $\Psi^{\rm
pendulum}_h$ is first order in~$h$. As in Section~\ref{sec:non_reso}, we
use $\Psi_{h}^{\rm pendulum}$ to build a symplectic and symmetric
scheme, denoted Algorithm~\ref{algo:pendulum_sym} in the sequel, following 
$$
(A,R,P_a,P_r) = \Psi_h^{\rm pend. sym.}(a,r,p_a,p_r) = 
\left( \Psi_{h/2}^{\rm pendulum} \right)^* \, \Psi_{h/2}^{\rm pendulum}
(a,r,p_a,p_r),
$$
where $\left( \Psi_h^{\rm pendulum} \right)^*$ denotes the the adjoint method,
which is denoted Algorithm~\ref{algo:pendulum_adjoint}.

\bookbox{
\begin{algorithm}[Symmetric Symplectic Scheme 
$\Psi^{\rm pend. sym.}_h(a,r,p_a,p_r)$]

Set $(a,r,p_a,p_r) = \left( a^n,r^n,p_a^n,p_r^n \right)$ and 
  perform the following steps:
\begin{enumerate}
\item Set $(\overline{A}, \overline{R}, \overline{P}_a,
  \overline{P}_r) = \Psi_{h/2}^{\rm pendulum} (a,r,p_a,p_r)$. 
\item Set $(A,R,P_a,P_r) = \left( \Psi_{h/2}^{\rm pendulum} \right)^*
(\overline{A}, \overline{R}, \overline{P}_a,\overline{P}_r)$.
\end{enumerate}
\quad Set
$
\left( a^{n+1},r^{n+1},p_a^{n+1},p_r^{n+1} \right) =
(A,R,P_a,P_r)  .
$
\label{algo:pendulum_sym}
\end{algorithm}
}

\bookbox{
\begin{algorithm}[Symplectic Scheme 
$\left( \Psi^{\rm pendulum}_h \right)^*(a,r,p_a,p_r)$]

Set $(a,r,p_a,p_r) = \left( a^n,r^n,p_a^n,p_r^n \right)$,
$\tau = h/\varepsilon$ and perform the following steps:
\begin{enumerate}
\item \text{Rotate the variables: set} 
$\dps
b = \frac{r}{\eps} \cos \tau +
p_r \sin \tau ,
\quad
p_b =
- \frac{r}{\eps} \sin \tau  + p_r \cos \tau .
$
\item Solve for $(A,B)$ in the equations
$$
\left\{
\begin{array}{rcl}
a &=& A - h p_a + 2 \eps^2 p_a (p_b \cos \tau + B \sin \tau) - 
2 \eps^2 p_b (p_a - h W'(A))
\\ \noalign{\vskip 4pt}
&-& \dps
h \varepsilon^2 \left[ \frac32 (B^2+p_b^2) (p_a - h W'(A)) - 
2 (p_a - h W'(A))^3 
\right],
\\ \noalign{\vskip 4pt}
b &=& \dps
B + \varepsilon p_a^2 \cos \tau - \eps (p_a - h W'(A))^2 
- \frac32 h \varepsilon p_b (p_a - h W'(A))^2 .
\end{array}
\right.
$$ 
\item Set
$\dps 
P_b 
= 
p_b 
+ \varepsilon p_a^2 \sin \tau 
- \frac32 h \varepsilon B (p_a - h W'(A))^2.
$
\item Set
\begin{eqnarray*}
P_a 
&=& 
\dps{
p_a - h W'(A) 
+ 2 h \eps^2 p_b (p_a - h W'(A)) W''(A) }
\\
&+& \dps{
h^2 \eps^2 \left[ \frac32 (B^2+p_b^2) (p_a - h W'(A)) - 2
  (p_a - h W'(A))^3 \right] W''(A). 
}
\end{eqnarray*}
\item Return to the original variables:
$
R = \eps B,
\quad
P_r = P_b.
$
\end{enumerate}
\quad Set
$
\left( a^{n+1},r^{n+1},p_a^{n+1},p_r^{n+1} \right) =
(A,R,P_a,P_r)  .
$
\label{algo:pendulum_adjoint}
\end{algorithm}
}

\subsection{Numerical results}
\label{sec:numerics_pendulum}

\begin{figure}[htbp]
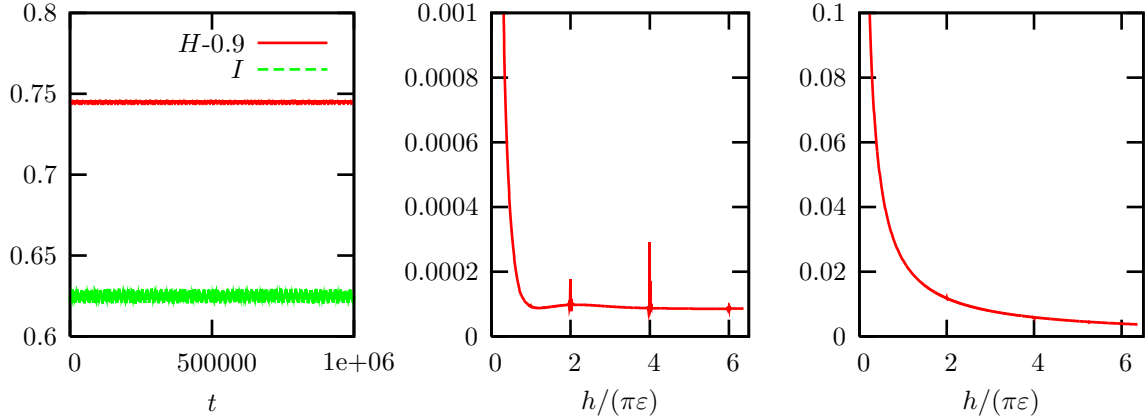

\centerline{
\input{figs/preservation_pendulum_small.tex}
\hspace{-7mm}
\input{figs/comparison_h0.02_pendulum_energy_small.tex}
\hspace{-7mm}
\input{figs/comparison_h0.02_pendulum_inv_adiab_small.tex}
}
\caption{\small
Left: energy (for convenience, we plot $H-0.9$) and adiabatic invariant $I$
along the trajectory computed with Algorithm~\ref{algo:pendulum_sym}
($\eps = 2.10^{-3}$ and $h=0.02$). 
Center and right: 
maximum variations (\ref{eq:estim}) of the energy (center) and 
of the adiabatic invariant $I$ (right) on the
time interval $[0,10^4]$, for several $\eps$ ($h = 0.02$), for Algorithm
\ref{algo:pendulum_sym}.} 
\label{fig:pendule}
\end{figure}

We have implemented Algorithm \ref{algo:pendulum_sym} on the
Hamiltonian~\eqref{eq:ham_pendule_int} with 
$$
W(a) = (\cos a)^2.
$$ 
We first choose $\eps=2.10^{-3}$ and $h=0.02$, and monitor the evolution
of the energy and of the adiabatic invariant 
\begin{equation}
\label{eq:I}
I = \frac{p_r^2}{2} + \frac{r^2}{2 \eps^2} 
\end{equation}
up to time $T=10^6$. Results
are shown on 
Fig.~\ref{fig:pendule}. We observe no drift.

\medskip

We next study the robustness of the algorithm as $\eps$ decreases.
 We set the time step to 
$h=0.02$, and consider the variations~\eqref{eq:estim} of the energy
and of the adiabatic invariant~\eqref{eq:I},  
over the time interval $t \in [0,10^4]$, for
stiffness $\eps$ varying between $10^{-3}$ to 1. Results are shown in
Fig.~\ref{fig:pendule}. 
Again, the algorithm performs equally
well when $\eps$ decreases to 0, up to extremely peaked resonances in the
energy preservation.

\section*{Acknowledgments}

This work was supported in part by the NSF Mathematical Sciences
Postdoctoral Research Fellowship, by the INRIA under the grant ``Action
de Recherche Collaborative'' HYBRID, and by the Agence Nationale de la
Recherche, under grant ANR-09-BLAN-0216-01 (MEGAS).


\begin{thebibliography}{10}

\bibitem{bornemann}
F.~Bornemann.
\newblock {\em Homogenization in time of singularly perturbed mechanical
  systems}.
\newblock Number 1687 in Lecture Notes in Mathematics. Springer, Berlin, 1998.

\bibitem{schuette}
F.~A. Bornemann and C.~Sch{\"u}tte.
\newblock Homogenization of {H}amiltonian systems with a strong constraining
  potential.
\newblock {\em Phys. D}, 102(1-2):57--77, 1997.

\bibitem{rennes}
F.~Castella, P.~Chartier, and E.~Faou.
\newblock An averaging technique for highly oscillatory {H}amiltonian problems.
\newblock {\em SIAM J. Numer. Anal.}, 47(4):2808--2837, 2009.

\bibitem{cohen_review}
D.~Cohen, T.~Jahnke, K.~Lorenz, and C.~Lubich.
\newblock Numerical integrators for highly oscillatory {H}amiltonian systems: a
  review.
\newblock In {\em Analysis, modeling and simulation of multiscale problems},
  pages 553--576. Springer, Berlin, 2006.

\bibitem{cras2}
M.~Dobson, C.~Le~Bris, and F.~Legoll.
\newblock Symplectic schemes for highly oscillatory {H}amiltonian systems with
  varying fast frequencies.
\newblock {\em C. R. Acad. Sci. Paris}, submitted.

\bibitem{feng}
K.~Feng.
\newblock Difference schemes for {H}amiltonian formalism and symplectic
  geometry.
\newblock {\em J. Comput. Math.}, 4(3):279--289, 1986.

\bibitem{mollify}
B.~Garc{\'{\i}}a-Archilla, J.~M. Sanz-Serna, and R.~D. Skeel.
\newblock Long-time-step methods for oscillatory differential equations.
\newblock {\em SIAM J. Sci. Comput.}, 20(3):930--963 (electronic), 1999.

\bibitem{grimm}
V.~Grimm and M.~Hochbruck.
\newblock Error analysis of exponential integrators for oscillatory
  second-order differential equations.
\newblock {\em J. Phys. A}, 39(19):5495--5507, 2006.

\bibitem{impulse1}
H.~Grubm\"uller, H.~Heller, A.~Windemuth, and K.~Schulten.
\newblock Generalized {V}erlet algorithm for efficient molecular dynamics
  simulations with long range interaction.
\newblock {\em Mol. Sim.}, 6(1-3):121--142, 1991.

\bibitem{HLW}
E.~Hairer, C.~Lubich, and G.~Wanner.
\newblock {\em Geometric numerical integration}, volume~31 of {\em Springer
  Series in Computational Mathematics}.
\newblock Springer-Verlag, Berlin, second edition, 2006.

\bibitem{cras1}
C.~Le~Bris and F.~Legoll.
\newblock Derivation of symplectic numerical schemes for highly oscillatory
  {H}amiltonian systems.
\newblock {\em C. R. Acad. Sci. Paris}, 344(4):277--282, 2007.

\bibitem{dcds}
C.~Le~Bris and F.~Legoll.
\newblock Integrators for highly oscillatory {H}amiltonian systems: an
  homogenization approach.
\newblock {\em Discrete Cont. Dyn-B}, 13(2):347--373, 2010.

\bibitem{mollify2}
J.~M. Sanz-Serna.
\newblock Mollified impulse methods for highly oscillatory differential
  equations.
\newblock {\em SIAM J. Numer. Anal.}, 46(2):1040--1059, 2008.

\bibitem{flavors}
M.~Tao, H.~Owhadi, and J.~E. Marsden.
\newblock Nonintrusive and structure preserving multiscale integration of stiff
  odes, sdes, and hamiltonian systems with hidden slow dynamics via flow
  averaging.
\newblock {\em Multiscale Modeling \& Simulation}, 8(4):1269--1324, 2010.

\bibitem{impulse2}
M.~Tuckerman, B.~Berne, and G.~Martyna.
\newblock Reversible multiple time scale molecular dynamics.
\newblock {\em J. Chem. Phys.}, 97(3):1990--2001, 1992.

\end{thebibliography}
\end{document}